\documentclass[11pt,a4paper]{amsart}
\usepackage[T1]{fontenc}
\usepackage{lmodern}
\usepackage{hyperref}
\usepackage{mathtools,bbm,amsthm,amsmath,amsfonts,amssymb,color,epic}
\usepackage{multirow,comment,todonotes}
\usepackage[export]{adjustbox}
\usepackage{tikz,pgfplots}
\usepackage{algorithmicx,algorithm,pifont}
\usepackage{microtype}
\usepackage{cite}
\pgfplotsset{compat=newest}
\usetikzlibrary{arrows,snakes,backgrounds,positioning}

\theoremstyle{plain}

\newcommand{\Sref}{S_0}
\newcommand{\Dref}{D_0}
\newcommand{\DefField}{{\bs\chi}}

\newcommand{\Rset}{\mathbb{R}}

\newcommand{\bs}[1]{{\boldsymbol#1}}
\newcommand{\xref}{\widehat{\bs x}}

\newcommand{\E}{\mathbb{E}}

\newcommand{\bfx}{{\bs x}}

\def\Letters{A,B,C,D,E,F,G,H,I,J,K,L,M,N,O,P,Q,R,S,T,U,V,W,X,Y,Z}
\makeatletter
\@for \@l:=\Letters \do{%
	\expandafter\edef\csname\@l bb\endcsname{\noexpand\ensuremath{\noexpand\mathbb{\@l}}}%
	\expandafter\edef\csname\@l bf\endcsname{{\noexpand\bf \@l}}%
	\expandafter\edef\csname\@l cal\endcsname{\noexpand\ensuremath{\noexpand\mathcal{\@l}}}%
	\expandafter\edef\csname\@l eu\endcsname{\noexpand\ensuremath{\noexpand\EuScript{\@l}}}%
	\expandafter\edef\csname\@l frak\endcsname{\noexpand\ensuremath{\noexpand\mathfrak{\@l}}}%
	\expandafter\edef\csname\@l rm\endcsname{{\noexpand\rm \@l}}%
	\expandafter\edef\csname\@l scr\endcsname{\noexpand\ensuremath{\noexpand\mathscr{\@l}}}%
}
\renewcommand{\d}{\operatorname{d\!}}
\newcommand{\Cor}{\operatorname{Cor}}
\newcommand{\Cov}{\operatorname{Cov}}

\newcommand{\isdef}{\mathrel{\mathrel{\mathop:}=}}

\definecolor{navy}{RGB}{102,153,255}
\definecolor{tuerkis}{RGB}{51,153,204}
\algrenewcommand\alglinenumber[1]{\ding{\numexpr191 + #1}}

\title[$p$-multilevel Monte Carlo for the acoustic scattering problem]
{$p$-multilevel Monte Carlo for acoustic scattering from
	large deviation rough random surfaces}

\author{J\"urgen D\"olz}
\email{doelz@ins.uni-bonn.de}
\author{Wei Huang}
\email{wei.huang@usi.ch}
\author{Michael Multerer}
\email{michael.multerer@usi.ch}

\address{
}

\begin{document}
	
\begin{abstract}
%%%%%%%%%%%%%%%%%%%%%%%%%%%%%%%%%%%%%%%%%%%%%%%%%%%%%%%%%%%%%%%%%%%%%%%%%%%%%%%%
We study time harmonic acoustic scattering on large deviation rough random 
scatterers. Therein, the roughness of the scatterers is caused by a low Sobolev 
regularity in the covariance function of their deformation field. The motivation 
for this study arises from physical phenomena where small-scale material defects 
can potentially introduce non-smooth deviations from a reference domain. The 
primary challenge in this scenario is that the scattered wave is also random, 
which makes computational predictions unreliable. Therefore, it is essential to 
quantify these uncertainties to ensure robust and well-informed design processes. 
While existing methods for uncertainty quantification typically rely on domain 
mapping or perturbation approaches, it turns out that large and rough random 
deviations are not satisfactory covered. To close this gap, and although counter 
intuitive at first, we show that the $p$-multilevel Monte Carlo
method can provide an efficient tool for uncertainty quantification in this 
setting. To this end, we discuss the stable implementation of higher-order 
polynomial approximation of the deformation field by means of barycentric 
interpolation and provide a cost-to-accuracy analysis. Our considerations are 
complemented by numerical experiments in three dimensions on a complex scattering 
geometry.
%%%%%%%%%%%%%%%%%%%%%%%%%%%%%%%%%%%%%%%%%%%%%%%%%%%%%%%%%%%%%%%%%%%%%%%%%%%%%%%%
\end{abstract}
	
%%%%%%%%%%%%%%%%%%%%%%%%%%%%%%%%%%%%%%%%%%%%%%%%%%%%%%%%%%%%%%%%%%%%%%%%%%%%%%%%
\maketitle
%%%%%%%%%%%%%%%%%%%%%%%%%%%%%%%%%%%%%%%%%%%%%%%%%%%%%%%%%%%%%%%%%%%%%%%%%%%%%%%%

\section{Introduction}
 
\subsection{Motivation}
We consider time harmonic, sound-soft acoustic scattering where the scatterer is 
a bounded obstacle. We assume that the scatterer $D(\omega)\subset\mathbb{R}^3$ 
is given as a large and rough random deviation of Karhunen-Lo\`eve-type from a 
reference domain $D_0$ with Lipschitz boundary. Here, we consider a deviation as 
rough if the covariance function of the deformation field has low Sobolev 
regularity. This setting is motivated by the fact that in reality defects in the 
scatterer's surface often happen on a relatively small scale and may even induce 
non-smooth artifacts in a previously smooth surface. Considering a plane incident 
wave $u_{\text{inc}}({\bs x}) 
= e^{-i\kappa\langle{\bs d},
	{\bs x}\rangle}$ with known wavenumber 
$\kappa$ and direction ${\bs d}$,
where $\|{\bs d}\|_2=1$, the total wave 
$
u = u_{\text{inc}}+u_{\text{s}}(\omega)
$
around the scatterer
is obtained by solving the exterior boundary value problem
\begin{equation}\label{eq:pde}
	%================================
	\begin{aligned}
		\Delta u + \kappa^2 u = 0\quad & \text{in}\ 
		\mathbb{R}^3\setminus\overline{D}(\omega), \\
		u = 0\quad                     & \text{on}\ S(\omega)=\partial D(\omega),\\
		\sqrt{r}\bigg(\frac{\partial u_{\mathrm{s}}(\omega)}{\partial r}-i\kappa
		u_{\mathrm{s}}(\omega)\bigg)
		\to 0\quad                     & \text{as}\ r = \|{\bs x}\|_2
		\to\infty,
	\end{aligned}
\end{equation}
to obtain the scattered wave $u_{\text{s}}(\omega)$. The sound-hard case can be 
treated in complete analogy, see for example \cite{DHM23}. 

The main difficulty in the described situation is that the scattered wave 
$u_{\text{s}}(\omega)$ is also random. To ensure robust and well-informed design 
processes in engineering and the applied sciences, it is of utmost importance to 
quantify these uncertainties. Due to the complexity of most underlying partial 
differential equations and the arising uncertainty, this needs to be achieved by 
means of numerical algorithms.

%%%%%%%%%%%%%%%%%%%%%%%%%%%%%%%%%%%%%%%%%%%%%%%%%%%%%%%%%%%%%%%%%%%%%%%%%%%%%%%%
\subsection{State-of-the-art}
Pioneering works on the acoustic scattering problem by deterministic scatterers 
with roughness on a small scale have been of interest for a long time with 
pioneering works going back to \cite{twersky1983reflection,kachoyan1987acoustic,
macaskill1988numerical}, see also \cite{darmon2020acoustic} and the references 
therein for a recent review. However, to the best of our knowledge, there are 
very scarce investigations into the acoustic scattering problem by rough random 
scatterers, especially in case of complex three-dimensional obstacles.

Within the mathematical community for uncertainty quantification, there have 
been significant contributions to develop algorithms quantifying this uncertainty 
in the solution to partial differential equations on random domains, mainly by 
computing statistical moments. These contributions can roughly be categorized into 
domain mapping approaches and perturbation approaches.
Domain mapping approaches \cite{HPS16,dolz2022isogeometric,DH2023,AJZ20,
JSZ17,HSSS2018,CNT2016,xiu} transfer the shape uncertainty into randomly varying 
coefficients on a fixed reference domain. This allows to deal with large domain 
deformations and the computation of statistical moments by evaluating the 
corresponding high-dimensional integrals, typically employing sampling based 
quadrature methods. Unfortunately, this requires the realization of a random 
domain and solving a partial differential equation thereon for each sample, 
which can be quite expensive. To reduce the number of required samples, 
quasi-Monte Carlo (QMC) methods \cite{Caf1998} and sparse grids \cite{BG2004a} 
exploit parametric smoothness and may be combined with multilevel sampling 
techniques \cite{Gil15,HPS12}.
Perturbation approaches compute the statistical moments by approximating them 
with truncated Taylor expansions in the Fr\'echet sense 
\cite{Dol2020,HSS2008a,DH2018,BN2014}. The additional correction terms can often 
be computed without the need for high-dimensional quadrature methods. Although 
perturbation approaches have the potential to deal with rough domain deformations, 
the obtained results are usually only reliable when the domain deformations are small.
A combination of the two approaches was presented in \cite{CNT2021,Mul18}, but relies 
on an a-prioriliy available splitting of the domain deformations into small and large 
amplitudes of the deformations. Summarizing, state-of-the-art methods rely on different 
kinds of smoothness and can either work with large and rather global or with small and 
rough randomness in the shape of the domain.

%%%%%%%%%%%%%%%%%%%%%%%%%%%%%%%%%%%%%%%%%%%%%%%%%%%%%%%%%%%%%%%%%%%%%%%%%%%%%%%%
\subsection{Challenges when dealing with rough random shapes}
The situation which is left out in the scope of currently available methods is 
the case when the domain deformations are rough and possibly large. This setting 
implies a number of challenges:
\begin{enumerate}
\item QMC and sparse grids do not work: The low Sobolev regularity leads to a 
slow decay of the eigenvalues of the associated covariance operator, leading to 
a true high-dimensionality of the parameter space. Indeed, the decay in the 
parametric dimension is such that QMC and sparse grid methods are no longer 
efficiently applicable.

\item Perturbation approaches are not suitable: The main assumption of 
perturbation approaches, that the amplitude of the perturbations is suitably 
small, is not necessarily fulfilled for rough domain deformations. The prime 
example are domain deformations described by the Mat\'ern-1/2 kernel, which is 
not differentiable on the diagonal but can model arbitrarily large domain 
deformations.
\item Resolving rough surfaces is expensive: It is well known that the numerical 
approximation of rough surface domains requires lower or lowest order polynomials 
on relatively refined meshes close to the roughness. Similar considerations hold 
for functions defined on these domains, and in particular for the solutions of 
partial differential equations thereon. This leads to a large number of degrees of 
freedom and thus to large computational cost per sample.
\end{enumerate}
In view of the above literature review, the only available option in this case is 
to use a Monte Carlo sampling strategy on a set of highly refined samples. However, 
this approach is computationally not feasible.

%%%%%%%%%%%%%%%%%%%%%%%%%%%%%%%%%%%%%%%%%%%%%%%%%%%%%%%%%%%%%%%%%%%%%%%%%%%%%%%%
\subsection{Contribution}
The purpose of this article is to illustrate on the example of acoustic 
scattering that the $p$-multilevel Monte Carlo method ($p$-MLMC), see 
\cite{BRFL+20}, has the potential to alleviate the aforementioned challenges. 
By $p$-MLMC we refer to a construction of the required multilevel hierarchy by 
increasing the polynomial degree of the ansatz spaces, rather than various mesh 
hierarchies. To this end, we
\begin{enumerate}
\item discuss the stable and efficient shape discretization by higher-order 
polynomials and barycentric interpolation,
\item perform a required sample to accuracy analysis for the $p$-MLMC,
\item provide numerical numerical examples in three-dimensions on a complex 
geometry.
\end{enumerate}
The use of a higher-order polynomial approximation method may seem unintuitive at 
first, but can be motivated quite naturally. To this end, we curiously observe that, 
on smooth surfaces, the eigenfunctions of covariance integral operators are always 
smooth, even if the covariance function itself is not. Thus, when the Karhunen-L\`eve 
expansion is truncated for numerical computations, the resulting deformation field, 
the deformed domain, and the solution to the scattering problem will always be smooth. 
We can thus expect exponential convergence in the polynomial degree if higher-order 
schemes are used.

%%%%%%%%%%%%%%%%%%%%%%%%%%%%%%%%%%%%%%%%%%%%%%%%%%%%%%%%%%%%%%%%%%%%%%%%%%%%%%%%
\subsection{Outline}
The remainder of this article is structured as follows. Section 
\ref{rough_random_domain_model} introduces the modeling of rough random domains 
using Mat\'ern kernels with low smoothness index and their representation by 
the Karhunen-Lo\`eve expansion. Section \ref{discretization_of_rough_random_domains} 
explores the conversion of a complex geometry, i.e., Michelangelo's David, 
from a triangular mesh to a NURBS representation and discusses 
its further high-order polynomial re-interpolation. Moreover,
the discretization of random deformation fields by barycentric interpolation is
introduced. In Section \ref{scattering_at_random_obstacles}, we revisit the boundary 
integral formulation of the problem and the computation of quantities of interest (QoI),
as described in \cite{dolz2022isogeometric}. 
Section \ref{p_Multilevel_Monte_Carlo_method} introduces the $p$-MLMC method and 
determines the optimal scaling factors for increasing 
polynomial degrees. Section \ref{numerical_experiments} presents extensive numerical 
studies conducted on a torus and the David geometry. Concluding remarks are stated
in Section~\ref{conclusion}.

%%%%%%%%%%%%%%%%%%%%%%%%%%%%%%%%%%%%%%%%%%%%%%%%%%%%%%%%%%%%%%%%%%%%%%%%%%%%%%%%
\section{Rough random domain model} \label{rough_random_domain_model}
%%%%%%%%%%%%%%%%%%%%%%%%%%%%%%%%%%%%%%%%%%%%%%%%%%%%%%%%%%%%%%%%%%%%%%%%%%%%%%%%
%=========================================================
\subsection{Modelling of random surfaces}
%=========================================================
Let \((\Omega,\Fcal,\Pbb)\) denote a complete probability space.
For the modelling of random domains, we assume the
existence of a Lipschitz continuous
reference surface \(\Sref\) bounding the reference domain \(\Dref\subset\Rbb^D\) 
and a random deformation field 
\[
	\DefField\colon\Omega \times \Sref \to \Rset^3
\]
such that there exists a constant \(C_{\operatorname{uni}} > 0\) with 
\begin{equation*}
	%==========================================
	\|\DefField(\omega)\|_{C^1(\Sref;\mathbb{R}^3)},
	\|\DefField^{-1}(\omega)\|_{%
		C^1(S(\omega);\mathbb{R}^3)}
	\leq C_{\operatorname{uni}} 
\end{equation*}
and 
\[
	S(\omega) = \DefField(\omega, \Sref) \quad
	\text{{for~\(\mathbb{P}\)-a.e.}\ }\omega\in\Omega.
\]
We make the practical assumption that the random deformation field $\DefField$ is
is known by its first- and second- order spatial statistics, i.e., 
its expected deformation field
\[
	\E[\DefField]\colon S_0\to\Rset^3
\]
and its matrix-valued covariance function
\[
	\Cov[\DefField]\colon S_0\times S_0\to\Rset^{3\times 3}.
\]
Without loss of generality, we assume that \(\E[\DefField]({\bs x})={\bs x}\). 
Otherwise, the reference domain has to be appropriately mapped, see \cite{HPS16}. 
The components of the covariance function is modeled by Mat\'ern covariance functions of 
some order $\nu > 0$. They are given by
\begin{equation*}
k_\nu(r) \isdef \frac{2^{1-\nu}}{\Gamma(\nu)}\bigg( \frac{\sqrt{2\nu}r}{l}\bigg)^\nu 
K_\nu\bigg( \frac{\sqrt{2\nu}r}{l}\bigg),
\end{equation*}
The eigenvalues of these kernels asymptotically decay like
\begin{equation}\label{eq:eig_dec}
	\lambda_j \leq Cj^{-(1+\frac{2\nu}{d})},
\end{equation}
where $d$ is the spatial dimension, see e.g.\ \cite{graham2015quasi}. 
In particular, the eigenvalues of the Gaussian kernel ($\nu=\infty$) decay exponentially. 
Later on, we shall employ the smooth Gaussian kernels to model smooth deformations, while
rougher kernels (with smaller $\nu$, and thus a slower eigenvalue decaying rate) are used 
to represent rough deformations. 

We remark that in our particular case of a two-dimensional manifold, when 
$\nu \leq 5$, neither the Halton sequence, see \cite{Caf98,Wan02}, nor the
anisotropic sparse grid quadrature, see \cite{HHPS18}, are dimension robust
due to the insufficient decay of the deformation field's Karhunen-Lo\`eve expansion.
Therefore, we exclusively consider the Monte Carlo method in this article.
%%%%%%%%%%%%%%%%%%%%%%%%%%%%%%%%%%%%%%%%%%%%%%%%%%%%%%%%%%%%%%%%%%%%%%%%%%%%%%%%
\subsection{Representation of random deformation fields}
%%%%%%%%%%%%%%%%%%%%%%%%%%%%%%%%%%%%%%%%%%%%%%%%%%%%%%%%%%%%%%%%%%%%%%%%%%%%%%%%
In this section, we recall the ideas presented in 
\cite{dolz2022isogeometric,huang2022isogeometric}.
As has been argued there, it is sufficient to compute
the Karhunen-Lo\`eve expansion of the random deformation
field exclusively with respect to the random surface and no
volume discretization is required at all.
Hence, given the expected deformation field \(\Ebb[\DefField]\)
and its matrix-valued covariance function \(\Cov[\DefField]\),
we can compute the surface Karhunen-Lo\`eve expansion
\begin{align*}
	\DefField(\omega,{\bs x})=\Ebb[\DefField]({\bs x})+       
	\sum_{k=1}^\infty\sqrt{\lambda_{k}}\DefField_k({\bs x}) 
	Y_{k}(\omega),\quad{\bs x}\in\Sref.                     
\end{align*}
Herein, the tuples \(\{(\lambda_{k}, \DefField_{k})\}_k\)
are the eigenpairs of the covariance operator
\begin{equation*}
	\begin{aligned}
		  & \mathcal{C}\colon [L^2(\Sref)]^3\to [L^2(\Sref)]^3, \\
		  & (\mathcal{C}{\bs v})({\bs x})\isdef\int_{\Sref}       
		\Cov[\DefField]({\bs x},{\bs x}')
		{\bs v}({\bs x}')\d\sigma_{{\bs x}'}
	\end{aligned}
\end{equation*}
and, for \(\lambda_k\neq 0\),
the centred and uncorrelated random variables
\(\{Y_k\}_k\) are given according to
\[
	Y_{k}(\omega)\isdef\frac{1}{\sqrt{\lambda_{k}}}\int_{\Sref}
	\big(\DefField({\omega},{\bs x})
	-\Ebb[\DefField]({\bs x})\big)^\intercal\DefField_{k}({\bs x})
	\d{\sigma_{{\bs x}}}.
\]
	
In practice, however, the random variables \(\{Y_k\}_k\)
are not known explicitly and need to be estimated.
We make the common model assumption that the random variables
\(\{Y_{k}\}_k\) are independent and uniformly distributed with
\(Y_{k}\sim\mathcal{U}(-1,1)\) for all \(k\). 
For numerical computations, the Karhunen-Lo\`eve expansion
has to be truncated after \(m\in\Nbb\) terms, where \(m\)
has to be chosen to meet an acceptable accuracy.
Then, by
identifying each random variable \(Y_k\) by its image
\(y_k\in[-1,1]\), we
arrive at the parametric deformation field
\begin{equation}\label{eq:parametricVectorFieldBoundary}
	\DefField({\bs y},{\bs x})=\Ebb[\DefField]({\bs x})+
	\sum_{k=1}^m\sqrt{\lambda_{k}}\DefField_{k}({\bs x})y_k,
	\quad{\bs y}\in\Gamma\isdef[-1,1]^{m}.
\end{equation}
The parametric deformation field gives rise to
the parametric surfaces
\begin{align}\label{eq:randomdomain}
	S({\bs y}) = \big\{\DefField({\bs y},{\bs x}): 
	{\bs x}\in \Sref\big\}.                        
\end{align}
	
%%%%%%%%%%%%%%%%%%%%%%%%%%%%%%%%%%%%%%%%%%%%%%%%%%%%%%%%%%%%%%%%%%%%%%%%%%%%%%%%
\section{Discretization of rough random domains} \label{discretization_of_rough_random_domains}
%%%%%%%%%%%%%%%%%%%%%%%%%%%%%%%%%%%%%%%%%%%%%%%%%%%%%%%%%%%%%%%%%%%%%%%%%%%%%%%%
%=========================================================
\subsection{Surface representation of complex geometries}\label{sec:surfrep}
%=========================================================
We assume the usual isogeometric setting for the surface
$S_0$ of the reference domain $D_0$. The surface $S_0$ can be decomposed 
into several smooth \emph{patches}
\[
	S_0 = \bigcup_{i=1}^M S_{0}^{(i)}.
\]
The intersection $S_{0}^{(i)}\cap S_{0}^{(i')}$ 
consists at most of a common vertex or a common edge for 
\(i\neq i^\prime\). In particular, each patch $S_{0}^{(i)}$ is
the image of an invertible NURBS mapping
\begin{equation}\label{eq:parametrization}
	{\bs s}_i\colon\square\to S_{0}^{(i)}
	\quad\text{with}\quad S_{0}^{(i)} = {\bs s}_i(\square)
	\quad\text{ for } i = 1,2,\ldots,M,
\end{equation}
where \({\bs s}_i\) is of the form
\begin{align*}
	\mathbf{s}_i(x,y)\isdef \sum_{0=i_1}^{k_1}\sum_{0=i_2}^{k_2}            
	\frac{\mathbf{c}_{i_1,i_2} b_{i_1}^{p_1}(x) b_{i_2}^{p_2}(y)            
	w_{i_1,i_2}}{ \sum_{j_1=0}^{k_1-1}\sum_{j_2=0}^{k_2-1} b_{j_1}^{p_1}(x) 
	b_{j_2}^{p_2}(y) w_{j_1,j_2}}                                           
\end{align*}
for control points $\mathbf{c}_{i_1,i_2}\in \Rset^3$ and weights 
$w_{i_1,i_2}>0$, where $ \lbrace b_{i_1}^{p_1} \rbrace_{0\leq i_1< k_1}$ and 
$ \lbrace b_{i_2}^{p_2} \rbrace_{0\leq i_2< k2}$ are the B-spline basis. 
We shall further follow the common convention that
parametrizations with a common edge coincide except for orientation.

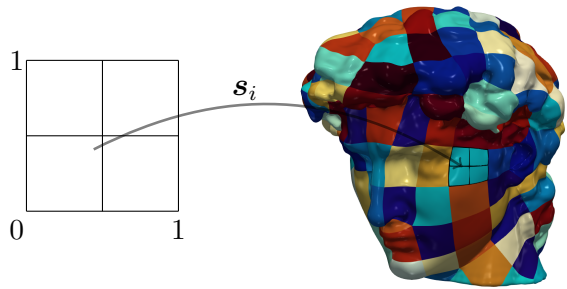
\begin{figure}[htb]
	\begin{center}
		\begin{tikzpicture}[
				scale=.25,
				axis/.style={thick, ->, >=stealth'},
				important line/.style={thick},
				every node/.style={color=black}
			]
			\draw (-.5,-1)node{{$0$}};
			\draw (8,-1)node{{$1$}};
			\draw (-.5,8)node{{$1$}};
			\draw (0,0)--(8,0);
			\draw (0,4)--(8,4);
			\draw (0,8)--(8,8);
			\draw (0,0)--(0,8);
			\draw (4,0)--(4,8);
			\draw (8,0)--(8,8);
			\draw (3,3)node(N1){};
			\draw (23.5,2)node(N2){};
			\draw
			(22.5,4)node{\includegraphics[width=0.5\textwidth,clip=true,trim=400 240
				280 190]{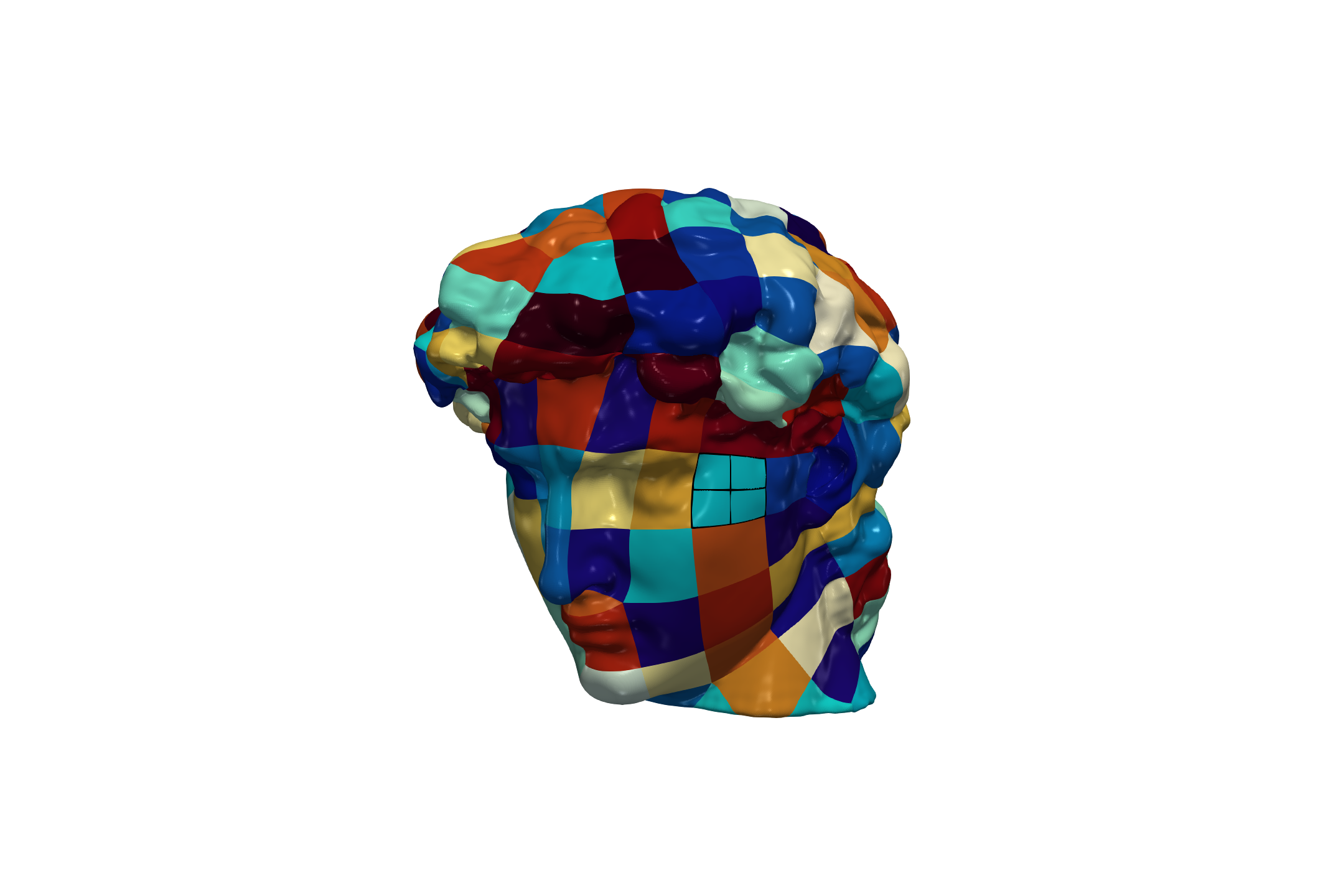}};
			\draw (11.5,6.4)node(N3){\large${\bs s}_i$};
			\path[->,line width=1pt, draw opacity=.5]  (N1) edge [bend left] (N2);
		\end{tikzpicture}
	\end{center}
	\caption{NURBS surface representation of Michelangelo's David. The spatial 
	refinement level is $\ell=1$ on \(\square\).}
	\label{fig:mesh}
\end{figure}

We remark that such s surface representation 
is often not directly available for more complex geometries. 
Instead, what we usually have are triangular meshes. 
The partitioning of surface meshes into conforming
quadrilateral patches can be solved by a workflow, 
which is based on determining a suitable set
of nodes called \emph{singularities} and connecting these by a set
of arcs called \emph{separatrices}, such that the resulting partition
consists of multiple valid conforming quadrilateral patches,
see \cite{dong2006spectral, campen2017partitioning} and the references therein. 
As a challenging geometry, we consider the head of Michelangelo's David are used 
in this paper. The initial geometry is a surface mesh from a 3D
scan (the scan is provided by the Statens Museum for Kunst under
the Creative Commons CC0 license).
Starting from a triangular mesh, we compute the critical points of the eigenfunctions 
of the Laplace-Beltrami operator. The critical points correspond to the geometric singularities. 
Next, we transform the triangular mesh into quadrangulations via the discrete 
Morse-Smale algorithm, which is available as \verb|Python| library 
\verb+Topology ToolKit+, see \cite{ttk,ttk19}. Finally, we obtain conforming quad patches by 
tracing separatrices starting from the critical points, followed by a subsequent 
approximation of each patch by NURBS surfaces, see Figure~\ref{fig:mesh}. Specifically, 
we use the 340-th largest eigenfunction of the Laplace-Beltrami operator on David's surface, 
which results in 187 patches in total. 
	
Regarding the boundary element solver, we define the ansatz space on \(S_0\) by patchwise 
lifting the standard B-spline basis space on the unit square. The same concept can be applied 
on the rough random surface \(S(\omega)\), see
\cite{huang2022isogeometric} for the details. The spatial level of uniform refinement 
on the unit square is denoted by $\ell$, and the degree of the B-spline used is indicated by $p$. 
Throughout the article, we fix $\ell=1$ and refine $p$ for increasing the accuracy of solutions. 
%%%%%%%%%%%%%%%%%%%%%%%%%%%%%%%%%%%%%%%%%%%%%%%%%%%%%%%%%%%%%%%%%%%%%%%%%%%%%%%%
\subsection{Barycentric interpolation} \label{sec:shape_discretization}
%%%%%%%%%%%%%%%%%%%%%%%%%%%%%%%%%%%%%%%%%%%%%%%%%%%%%%%%%%%%%%%%%%%%%%%%%%%%%%%%
We focus on the situation that the NURBS patches under
consideration are smooth, i.e., polynomial.
Therefore, we may employ a polynomial model
to represent the random deformation field.
We start by re-interpolating 
the patches \({\bs s}_i\colon\square\to\Sref\).
To this end, given a polynomial degree \(q\),
we introduce the Chebyshev points
of the second kind 
\begin{equation}\label{eq:chebNods}
	\xi_k\isdef\frac 1 2\bigg( 1- \cos\frac{k\pi}{q}\bigg),\quad k=0,\ldots, q.
\end{equation}
We denote the corresponding Lagrange polynomials by
\[
	L_k(s)\isdef\prod_{j\neq k}\frac{s-\xi_j}{\xi_k-\xi_j}.    
\]
The corresponding interpolation points on the unit square are then
given by a tensor product construction according to
\[
	{\bs\xi}_{k,k'}\isdef\begin{bmatrix}\xi_k\\ \xi_k'\end{bmatrix}\in\square.
\]
The mapped points
\[
	{\bs p}_{i,k,k'}\isdef{\bs s}_i({\bs\xi}_{k,k'})
\]
serve as landmark points for the representation of the
random surface. The set of all landmark points is denoted by
\begin{equation}\label{eq:landmarks}
{\bs\Theta}\isdef\{{\bs p}_{i,k,k'}: 1\leq i\leq M, 0\leq k,k'\leq q\}.
\end{equation}

The deformation field is now represented according to
\[
	\DefField\big(\omega,{\bs s}_i(s,t)\big)
	\approx\sum_{k=0}^{q}\sum_{k'=0}^{q}
	\DefField(\omega,{\bs p}_{i,k,k'})L_{k}(s)L_{k'}(t).
\]
To obtain a stable polynomial interpolation also for a large
number of landmark points, we employ the barycentric interpolation
formula from \cite{berrut2004barycentric}. Factoring out the node polynomial
and regarding that the Lagrange polynomials form a partition of unity
yields the one-dimensional polynomial interpolant of the function
\(f\colon[0,1]\to\Rbb\) by means of the representation
\begin{equation}\label{eq:uniInterp}
	p(s)=\rho(s)\sum_{k=0}^q\frac{w_k}{s-\xi_k}f(\xi_k),\quad
	\rho(s)\isdef\bigg(\sum_{k=0}^q\frac{w_k}{s-\xi_k}\bigg)^{-1}.
\end{equation}
Herein, the weights are given by
\[
	w_k\isdef\prod_{j\neq k}\frac{1}{\xi_k-\xi_j}.
\]
For the Chebyshev points of the second kind, the weights
can be precomputed by the formula, cp.\ \cite{salzer1972lagrangian},
\[
	w_k=(-1)^k\delta_k,\quad\text{where}\ 
	\delta_k\isdef\begin{cases}0.5,&\text{if }k=0,q,\\ 
	1,& \text{else}.
	\end{cases}
\]

Obviously, the weights are
invariant under affine transforms, as any such transforms
cancels out due to the barycentric formula. Moreover, in
the concrete case, the weights are also invariant under
the reversal of the order of the interpolation nodes,
as we have tacitly done in the definition \eqref{eq:chebNods}.
Now, tensorizing the univariate interpolant \eqref{eq:uniInterp}
yields the final representation
\begin{equation}\label{eq:reinterpolation}
	\DefField\big(\omega,{\bs s}_i(s,t)\big)
	\approx
	\rho(s)\rho(t)\sum_{k,k'=0}^q\frac{w_k}{s-\xi_k}
	\frac{w_k'}{s-\xi_{k'}}
	\DefField(\omega,{\bs p}_{i,k,k'}).
\end{equation}

\begin{figure}
	\begin{center}
		\begin{tikzpicture}[
				scale=.2,
				axis/.style={thick, ->, >=stealth'},
				important line/.style={thick},
				every node/.style={color=black}
			]
			\draw(3.5,3.6)node{%
				\includegraphics[scale=0.038,clip=true,trim=400 200 400 200]{%
			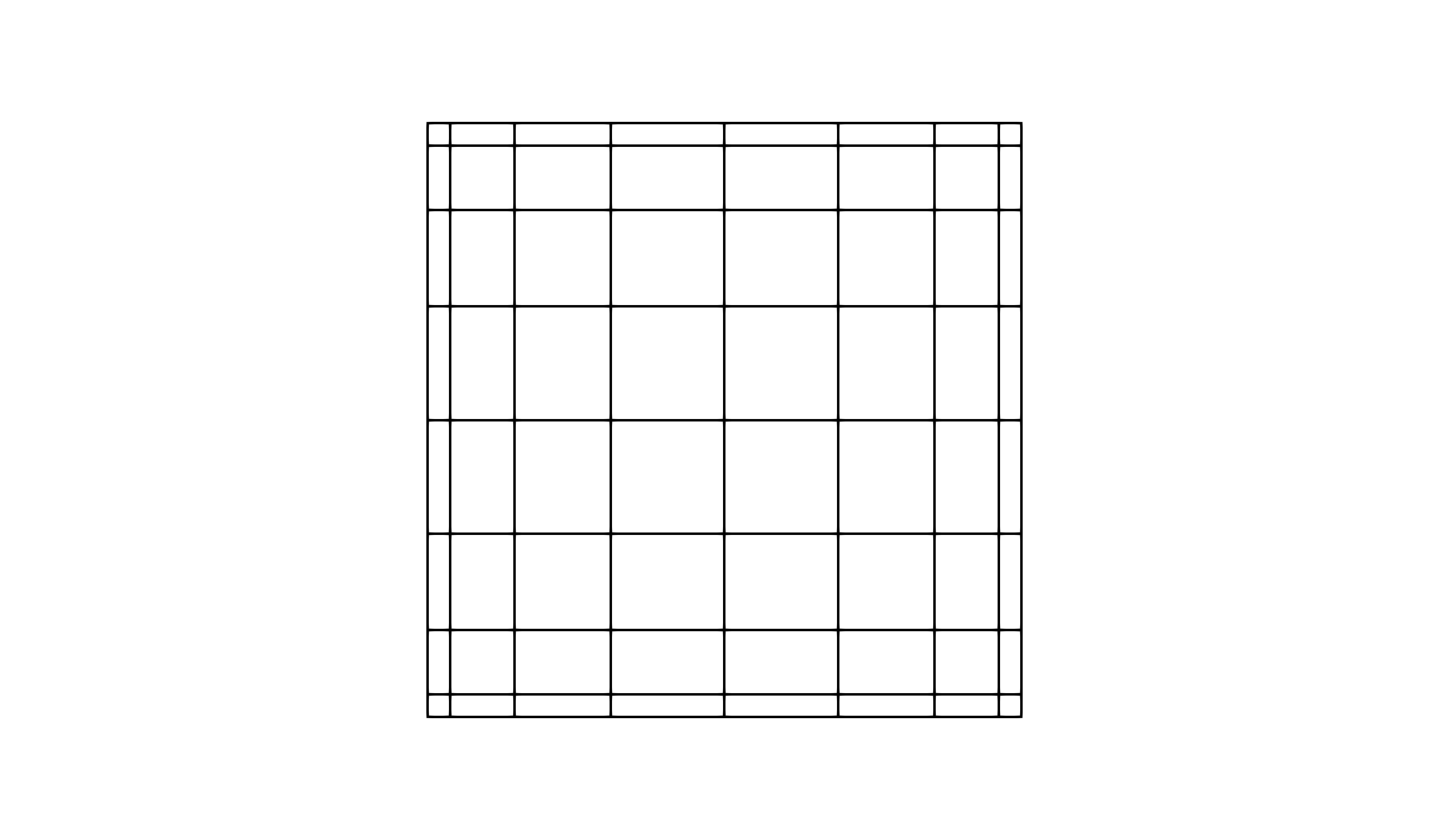}};
			\draw (-1,-1.5)node{{$0$}};
			\draw (7.5,-1.5)node{{$1$}};
			\draw (-1,8)node{{$1$}};
			\draw (3,3)node(N1){};
			\draw (28,1.5)node(N2){};
			\draw(25.5,3.6)node{%
				\includegraphics[scale=0.12,clip=true,trim=500 200 520 200]{%
			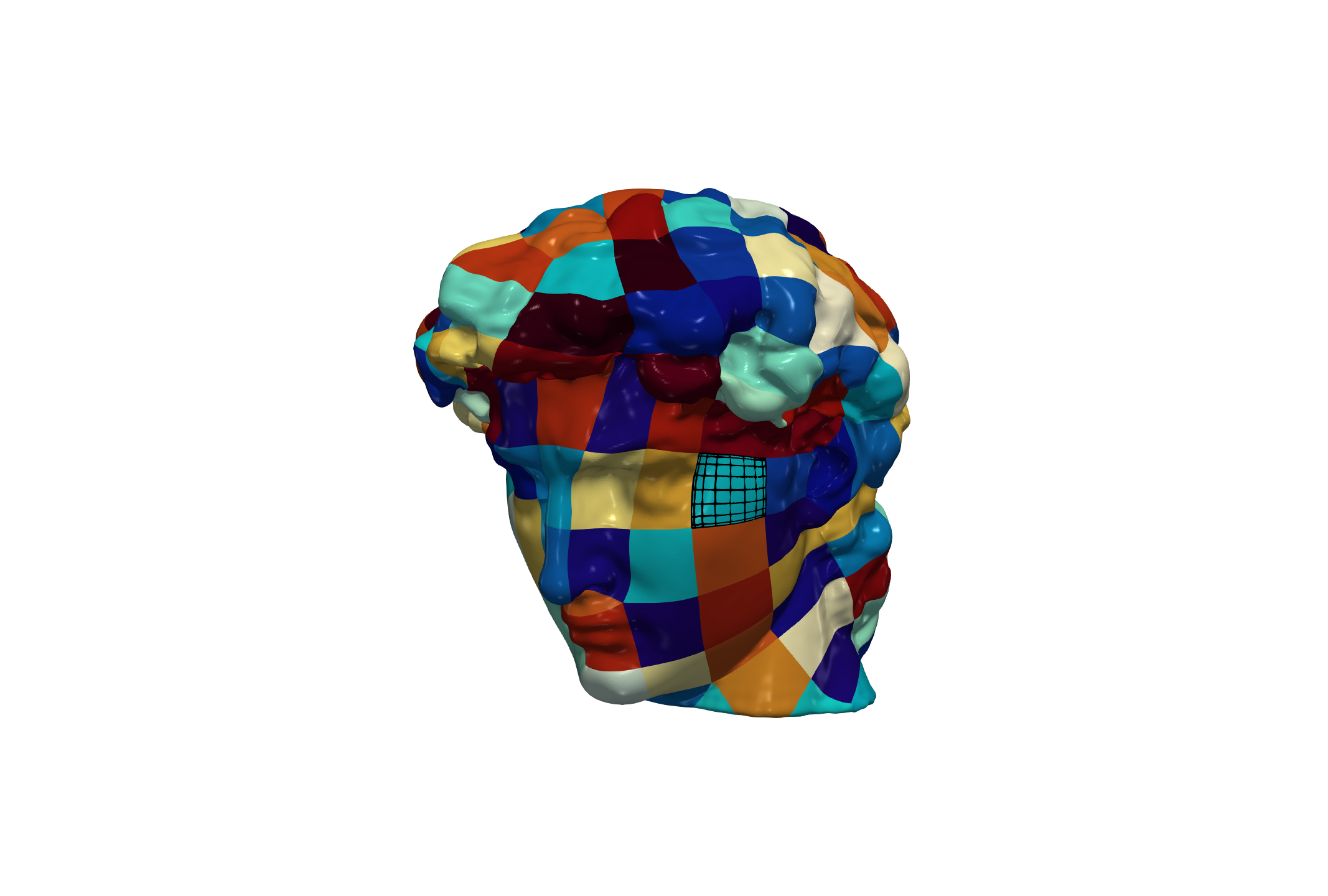}};
			\draw (14.7,6.8)node(N3){${\bs s}_i$};
			\path[->,line width=1pt]  (N1) edge [bend left] (N2);
		\end{tikzpicture}
		\vspace*{-1.5em}
	\end{center}
	\caption{\label{fig:p-reinterpolation}Tensor product of the 
	Chebyshev points of the second kind and surface re-interpolation 
	by polynomials of degree $q=8$.}
\end{figure}
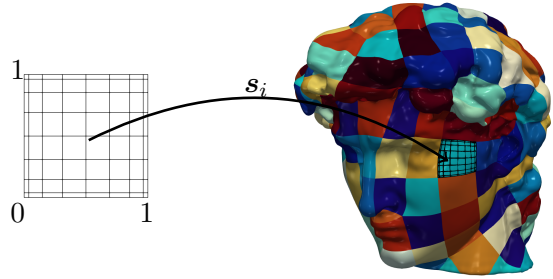

Here, the polynomial degree \(q\) is always chosen
such that the overall consistency error is met. 
An example of the re-interpolation of David is shown in 
Figure~\ref{fig:p-reinterpolation}.
	
%%%%%%%%%%%%%%%%%%%%%%%%%%%%%%%%%%%%%%%%%%%%%%%%%%%%%%%%%%%%%%%%%%%%%%%%%%%%%%%%
\subsection{Discretization of random deformation fields} \label{sec:field_discretization}
%%%%%%%%%%%%%%%%%%%%%%%%%%%%%%%%%%%%%%%%%%%%%%%%%%%%%%%%%%%%%%%%%%%%%%%%%%%%%%%%
The landmark points for the re-interpolation are also used to discretize the covariance operator. 
To evaluate the Karhunen-Lo\`eve expansion at the landmark points in \({\bs\Theta}\)
from \eqref{eq:landmarks}, 
we have to solve the eigenvalue problem for the covariance matrix
\[
	{\bs C}\isdef
	\big[\operatorname{Cov}[{\bs\chi}](\bfx,\bfx')
	\big]_{\bfx,\bfx'\in{\bs\Theta}}
	\in\Rbb^{3n\times 3n},
\]
where we set \(n\isdef|{\bs\Theta}|\).
Because solving the eigenvalue problem is of cost \(\Ocal(n^3)\), we first use 
a low-rank approach via the pivoted Cholesky decomposition, see \cite{HPS}. 
Having \({\bs C}\approx{\bs L}{\bs L}^\intercal\),
\(\operatorname{rank}{\bs L}=m\ll 3n\), we solve
\[
{\bs L}{\bs L}^\intercal{\bs v}=\lambda{\bs v}\quad\text{via}\quad
{\bs L}^\intercal{\bs L}=\tilde{\bs V}{\bs\Lambda}\tilde{\bs V}^\intercal.\]
Given that \((\lambda_,\tilde{\bs v}_i)\) is an eigenpair of the latter,
\((\lambda_k,{\bs L}\tilde{\bs v}_k)\)
is an eigenpair of the former. Finally, we
obtain the representation of the random deformation at the landmarks defined by 
\begin{equation}
	\DefField(\omega,{\bs\Theta})=
	{\bs\Theta}+\alpha\cdot\texttt{reshape}\big({\bs L}\tilde{\bs V}{\bs y}, 3, n\big).
	\label{eq:update}
\end{equation}
Herein, $\alpha$ is the magnitude of the random deformation, which steers the impact 
of the random deformation, and $\bs y$ satisfies the uniform distribution in \([-1,1]^m\). 

In summary, having a NURBS representation of the scatterer's surface, we first 
re-interpolate each patch at once with a tensor product of Lagrange polynomials of degree $q$. 
Then, we compute a low-rank approximation of of the rough random deformation field evaluated
at the landmark points via the pivoted Cholesky decomposition. Finally, the geometry deformations 
are obtained by updating the landmark points using \eqref{eq:update} and re-interpolation via 
\eqref{eq:reinterpolation}.
	
%%%%%%%%%%%%%%%%%%%%%%%%%%%%%%%%%%%%%%%%%%%%%%%%%%%%%%%%%%%%%%%%%%%%%%%%%%%%%%%%
\section{Scattering at random obstacles} \label{scattering_at_random_obstacles}
%%%%%%%%%%%%%%%%%%%%%%%%%%%%%%%%%%%%%%%%%%%%%%%%%%%%%%%%%%%%%%%%%%%%%%%%%%%%%%%%
%=========================================================
\subsection{Boundary integral equations}
%=========================================================
In this section, we recall the computation of scattered waves
\eqref{eq:pde} for each realization of the rough random scatterer 
using boundary integral equations. Throughout the paper, 
the scattered wave $u_s(\omega)$, the total wave $u(\omega)$, 
the scatterer $D(\omega)$, and the corresponding surface $S(\omega)$ 
are all dependent on the stochastic term $\omega$. 
For simplicity of notation, we omit the term $\omega$ in this paragraph.
	
We introduce the acoustic single layer operator
\[
	\mathcal{V}\colon H^{-1/2}(S)\to H^{1/2}(S),\quad
	(\mathcal{V}\rho)({\bs x})\isdef\int_S
	\Phi({\bs x},{\bs z})\rho({\bs z})\d\sigma_{\bs z},
\]
and the acoustic double layer operator
\[
	\mathcal{K}\colon L^2(S)\to L^2(S),\quad
	(\mathcal{K}\rho)({\bs x})\isdef\int_S
	\frac{\partial\Phi({\bs x},{\bs z})}
	{\partial{\bs n}_{\bs z}}\rho({\bs z})\d\sigma_{\bs z}.
\]
Herein, ${\bs n}_{\bs z}$ denotes the outward pointing normal vector
at the surface point ${\bs z}\in S$,  while $\Phi(\cdot,\cdot)$ denotes 
the fundamental solution for the Helmholtz equation. In our case, 
the fundamental function is given by
\[
	\Phi({\bs x}, {\bs z})
	= \frac{e^{-i\kappa\|{\bs x}-{\bs z}\|_2}}{4\pi\|{\bs x}-{\bs z}\|_2}.
\]
Considering an incident plane wave $u_{\mathrm{inc}}({\bs x}) = e^{-i\kappa\langle{\bs d},
	{\bs x}\rangle}$, $\|{\bs d}\|_2=1$, the Neumann data of the total 
	wave $u=u_{\mathrm{inc}}+u_{\mathrm{s}}$ at the 
surface $S$ can be determined by the boundary integral equation
\begin{equation*}
	\left(\frac{1}{2} + \mathcal{K^\star} - i\eta\mathcal{V}\right)
	\frac{\partial u}{\partial {\bs n}} 
	= \frac{\partial u_{\mathrm{inc}}}{\partial {\bs n}} - i\eta u_{\mathrm{inc}}
	\quad\text{on $S$},
\end{equation*}
with $\eta=\kappa/2$, cp.\ \cite{CK2}. 
	
\begin{figure}[ht]
	\begin{tikzpicture}[
			scale=.25,
			axis/.style={thick, ->, >=stealth'},
			important line/.style={thick},
			every node/.style={color=black}
		]
		\draw (0,0)node(N0){};
		\draw (8,0)node(N1){};
		\path[->,line width=1pt]  (N0) edge (N1);
		\draw (4,2)node{$u_{\mathrm{inc}}$};
		\draw (20,0)node{\includegraphics[scale=.14, clip, trim= 600 320 300 320]{%
		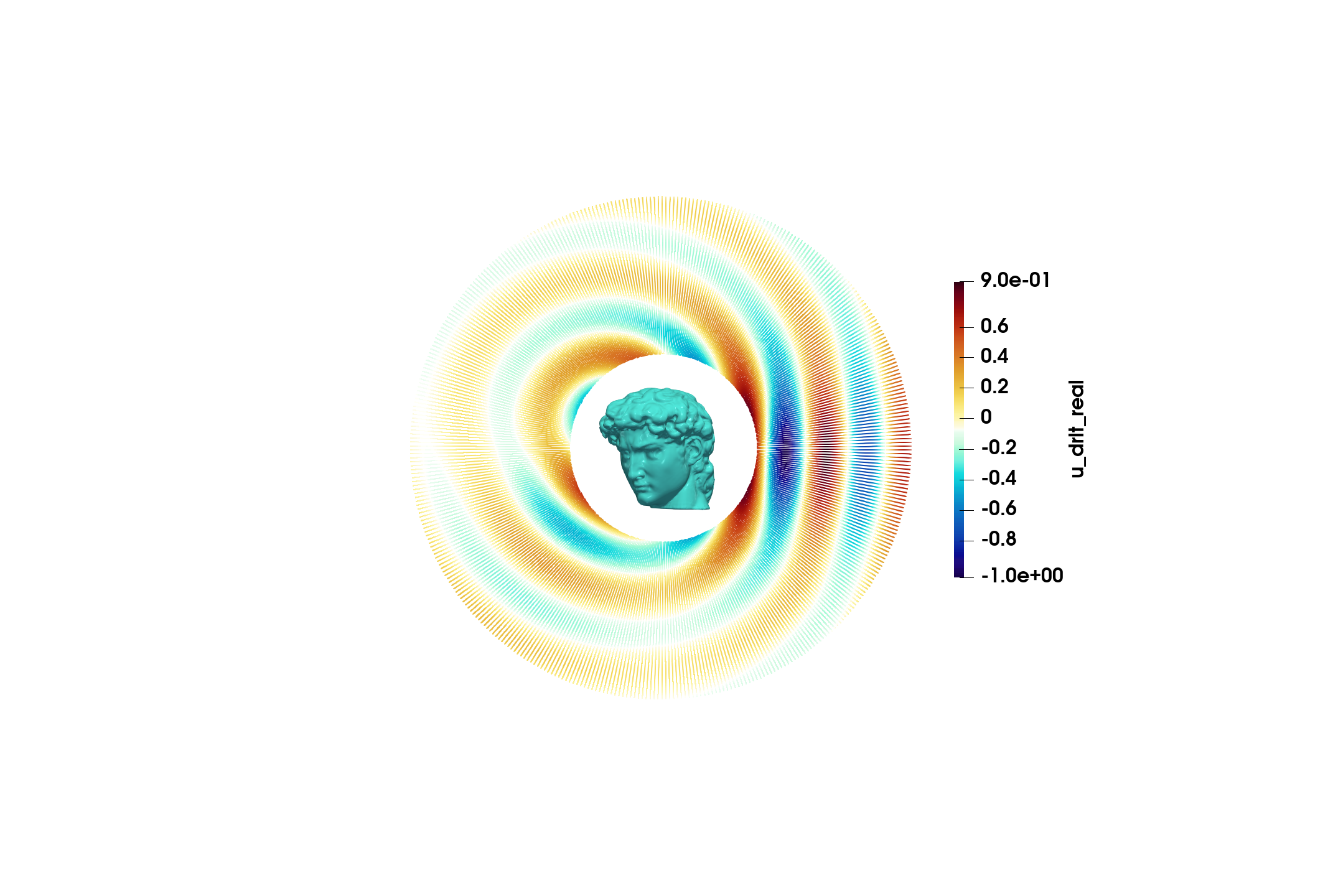}};
	\end{tikzpicture}
	\caption{It shows the real part of the scattering waves by Davide in the 
	case of an incident plane wave with wavenumber $\kappa=5$ and 
	direction ${\bs d} = [1,0,0]^\intercal$.}
	\label{fig:interface}
\end{figure}
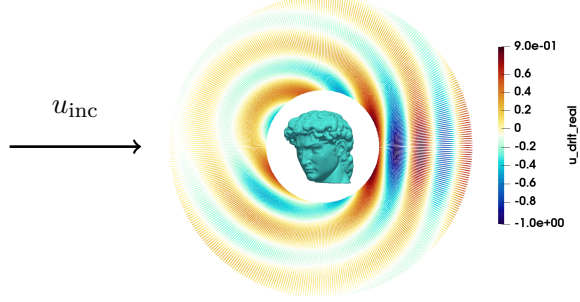

From the Cauchy data of $u$ at $S$, we 
can determine the scattered wave $u_{\mathrm{s}}$ at any exterior point 
outside the obstacle by applying the potential evaluation
\begin{equation}\label{eq:solution1}
	u_{\mathrm{s}}({\bs x}) = \int_S\Phi({\bs x}, {\bs z})
	\frac{\partial u}{\partial{\bs n}_{\bs z}}({\bs z})\d\sigma_{\bs z},
	\quad {\bs x}\in\mathbb{R}^3\setminus\overline{D}.
\end{equation}
In Figure \ref{fig:interface}, we visualise one example of a scattered wave by Michelangelo's David.	
	
%%%%%%%%%%%%%%%%%%%%%%%%%%%%%%%%%%%%%%%%%%%%%%%%%%%%%%%%%%%%%%%%%%%%%%%%%%%%%%%%
\subsection{Scattered wave representation at an artificial interface}
%%%%%%%%%%%%%%%%%%%%%%%%%%%%%%%%%%%%%%%%%%%%%%%%%%%%%%%%%%%%%%%%%%%%%%%%%%%%%%%%
As proposed in \cite{dolz2022isogeometric}, we introduce an artificial interface $T\subset\mathbb{R}^3$, 
being sufficiently large to guarantee that $T$ encloses all realizations of
the domain $D$. Applying \eqref{eq:solution1} yields the Dirichlet data of the scattered 
wave at the artificial interface, i.e., $u_{\mathrm{s}}|_T$.
With regard to the Neumann data $(\partial u_{\mathrm{s}}/\partial {\bs n})|_T$, 
we may take the normal derivative of \eqref{eq:solution1} at $T$, which yields
\[
	\frac{\partial u_{\mathrm{s}}}{\partial {\bs n}_{\bs x}}({\bs x}) 
	= \int_S \frac{\partial\Phi({\bs x}, {\bs z})}
	{\partial {\bs n}_{\bs x}}\frac{\partial u}
	{\partial{\bs n}_{\bs z}}({\bs z})\d\sigma_{\bs z},
	\quad {\bs x}\in T.
\]
	
Having the Cauchy data on $T$, the representation formula for $T$ is given by
\begin{equation}\label{eq:solution2}
	%==============================================
	u_{\mathrm{s}}({\bs x}) = \int_T \bigg\{
	\frac{\partial\Phi({\bs x}, {\bs z})}{\partial {\bs n}_{\bs z}}
	u_{\mathrm{s}}({\bs z})
	-
	\Phi({\bs x}, {\bs z})
	\frac{\partial u_{\mathrm{s}}}{\partial{\bs n}_{\bs z}}({\bs z})
	\bigg\}\d\sigma_{\bs z},
\end{equation}
see \cite{CK2}.
Then, we may evaluate scattered wave outside $T$ applying \eqref{eq:solution2} for 
$T$ instead of \eqref{eq:solution1} for $S$.
	
The significant advantage of \eqref{eq:solution2} over \eqref{eq:solution1} is 
that the artificial interface and the function space at $T$  are fixed in contrast 
to the shape of the random obstacle and non-nested hierarchical B-spline function 
spaces of increasing polynomial degrees, respectively. 
	
%%%%%%%%%%%%%%%%%%%%%%%%%%%%%%%%%%%%%%%%%%%%%%%%%%%%%%%%%%%%%%%%%%%%%%%%%%%%%%%%
\subsection{Quantities of interest}
%%%%%%%%%%%%%%%%%%%%%%%%%%%%%%%%%%%%%%%%%%%%%%%%%%%%%%%%%%%%%%%%%%%%%%%%%%%%%%%%
It is sufficient to compute the first and second moments of the 
Cauchy data at the artificial interface $T$ for 
knowing the statistics at the points outside $T$. We recall the formulas here.
For any ${\bs x}\in\mathbb{R}^3$ lying outside the interface 
$T$, it holds
\begin{equation*}
	\begin{aligned}
		\E[u_{\mathrm{s}}]({\bs x}) = \int_T \bigg\{
		  & \frac{\partial\Phi({\bs x}, {\bs z})}{\partial {\bs n}_{\bs z}} 
		\E[u_{\mathrm{s}}]({\bs z})
		  -\Phi({\bs x}, {\bs z})                                         
		\E\bigg[\frac{\partial u_{\mathrm{s}}}{\partial{\bs n}_{\bs z}}\bigg]({\bs z})\bigg\}
		\d\sigma_{\bs z}.
	\end{aligned}
\end{equation*}
With respect to the covariance, we find for two points ${\bs x},{\bs x}'\in\mathbb{R}^3$ lying 
outside of the interface $T$ the deterministic expression 
\begin{equation*}
	\begin{aligned}
		\Cor[u_{\mathrm{s}}]({\bs x},{\bs x}') & = \int_T\int_T                                                              
		\bigg\{\Phi({\bs x}, {\bs z})\overline{\Phi({\bs x}', {\bs z}')}
		\Cor\!\!\bigg[
		\frac{\partial u_{\mathrm{s}}}{\partial{\bs n}}\bigg]({\bs z},{\bs z}')\\
		                                       &\phantom{=}\qquad\qquad - \Phi({\bs x}, {\bs z})\overline{\frac{\partial\Phi({\bs x}', {\bs z}')} 
		{\partial {\bs n}_{{\bs z}'}}}
		\Cor\!\!\bigg[
		\frac{\partial u_{\mathrm{s}}}
		{\partial{\bs n}},u_{\mathrm{s}}\bigg]({\bs z},{\bs z}')\\
		                                       &\phantom{=}\qquad\qquad - \frac{\partial\Phi({\bs x}, {\bs z})}                                   
		{\partial {\bs n}_{\bs z}}\overline{\Phi({\bs x}', {\bs z}')}
		\Cor\!\!\bigg[u_{\mathrm{s}},\frac{\partial u_{\mathrm{s}}}
		{\partial{\bs n}}\bigg]({\bs z},{\bs z}')\\
		                                    &\phantom{=}\qquad\qquad+ 
		                                    \frac{\partial\Phi({\bs x}, {\bs z})}{\partial {\bs n}_{\bs z}}             
		\overline{\frac{\partial\Phi({\bs x}', {\bs z}')}
		{\partial {\bs n}_{{\bs z}'}}}
		\Cor[u_{\mathrm{s}}]({\bs z},{\bs z}')\bigg\}
		\d\sigma_{{\bs z}'}\d\sigma_{\bs z}.
	\end{aligned}
\end{equation*}
	
%%%%%%%%%%%%%%%%%%%%%%%%%%%%%%%%%%%%%%%%%%%%%%%%%%%%%%%%%%%%%%%%%%%%%%%%%%%%%%%%
\section{$p$-Multilevel Monte Carlo method} \label{p_Multilevel_Monte_Carlo_method}
%%%%%%%%%%%%%%%%%%%%%%%%%%%%%%%%%%%%%%%%%%%%%%%%%%%%%%%%%%%%%%%%%%%%%%%%%%%%%%%%
To compute QoIs, we rely on the Monte Carlo method, whose 
expected error 
is proportional to $N^{-1/2}$, where $N$ is the number of samples,
given that the integrand has bounded variance. 
To speed up the Monte Carlo method, we employ the 
$p$-MLMC method. This approach constructs a sequence of control
variates for variance reduction of the integrand. This variance reduction is
achieved by the decreasing spatial approximation error with respect to the polynomial
degree $p$. This error decreases even exponentially given that all involved
quantities are smooth.
	
The sequence of control variates in the multilevel Monte Carlo method is 
based on a telescoping sum, see e.g.\ \cite{H2, giles2008multilevel,Gil15}. 
We adapt their approach to the multilevel method for $p$-refinement, which yields the and
approximation for each of the two QoIs. In both cases we refer to
\(\mathcal{Q}_{P-p}\) to the Monte Carlo quadrature operator associated to the
polynomial degree $p$, where \(P\) is the maximum degree.
For the expectation of the quantity \(\rho^{(P)}\), we employ the estimator
	      \[
	      	\mathbb{E}[\rho]({\bs z})
	      	\approx\sum_{p=0}^P
	      	\mathcal{Q}_{P-p}\big(\rho^{(p)}(\cdot,{\bs z})-
	      	\rho^{(p-1)}(\cdot,{\bs z})\big),
	      \]
while its correlation is computed by
\[
	      	\Cor[\rho\otimes\mu]({\bs z},{\bs z}')  \approx
	      	\sum_{p=0}^P
	      	\mathcal{Q}_{P-p}\big((\rho\otimes\mu)^{(p)}(\cdot,{\bs z},{\bs z}')
	      	-(\rho\otimes\mu)^{(p-1)}(\cdot,{\bs z},{\bs z}')\big).
	      \]

The sequence of $\{\rho^{(p)}\}$ approximates $\rho$ with 
not only increasing accuracy but also increasing cost. Besides, 
we assume the variance of the difference $\rho^{(p)} - \rho^{(p-1)}$ 
between two consecutive levels decreases with respect to $p$. 
We recall that in $p$-refinement, the levels indicate the degrees of B-spline basis. 
By computing relatively few samples at the high levels but lots of samples at the coarse levels, 
we substantially save in terms of the total computations, see \cite{Gil15}. 
In particular, if we use $N_p$ samples for \(\mathcal{Q}_{P-p}\), we can rewrite 
the formula for the expectation according to
\begin{equation*}
	\mathbb{E}[\rho]({\bs z}) \approx N_0^{-1}\sum_{n=1}^{N_0} 
	\rho^{(0)}({\bs z},\omega_{0,n}) + \sum_{p=1}^P N_p^{-1}
	\sum_{n=1}^{N_p}\big(\rho^{(p)}({\bs z},\omega_{p,n}) - \rho^{(p-1)}({\bs z},\omega_{p,n})\big).
\end{equation*}
The same concept can be applied to the computation of the correlation. 
Let $c_0, v_0$ denote the cost and the variance of $\rho^{(0)} $, respectively. 
For $p>0$, $c_p, v_p$ are the cost and variance of $\rho^{(p)} - \rho^{(p-1)}$, accordingly. 
Then, the overall cost and variance of multilevel Monte Carlo estimate are 
\[
C(P)\isdef\sum_{p=0}^P N_p c_p
\] and 
\[
V(P)\isdef\sum_{p=0}^PN_{p}^{-1}v_p,
\] respectively, see \cite{Gil15}. Finding the minimum value of the 
overall cost subject to the constraint that the overall variance is $V(P)=\epsilon^2$ 
at most is equivalent to finding the values of $\{N_p\}$ and Lagrange multipliers 
$\lambda$ that simultaneously satisfy the equations defined by
	
\begin{equation}\label{eq:opt}
	\begin{cases}
		\nabla_{N_p}C(P) = \lambda \nabla_{N_p}\big(V(P) - \epsilon^2\big) \quad \text{for} \quad p=0,\cdots,P, \\
		                                                                                                                                    \\
		V(P) = \epsilon^2.                                                                                               
	\end{cases}       
\end{equation}
By solving the linear system of equations \eqref{eq:opt}, we obtain that when 
\begin{equation}\label{eq:opt_sol}
N_p \sim \sqrt{\frac{\lambda v_p}{c_p}}, \quad \text{where} \quad \lambda = \left(\frac{\sum_{p=0}^P \sqrt{c_pv_p}}{\epsilon^2}\right)^2,
\end{equation}
the minimum overall cost is achieved, see \cite{Gil15}. The computation at the low levels is cheaper
compared to calculations at the high levels, i.e., $c_p$ is a monotonically increasing function of the 
polynomial degree $p$, and the variance $v_p$ is a monotonically decreasing function of degree $p$. 
Equation \eqref{eq:opt_sol} suggests that the optimal number of samples decreases like 
$\sqrt{v_p/c_p}$ with increasing degrees. Specifically, when increasing from 
$p-1$ to $p$, the number of samples should decrease by a factor of 
\[
	\gamma_p \isdef \sqrt{\frac{v_{p-1}}{v_p}}\sqrt{\frac{c_p}{c_{p-1}}}.
\]
%%%%%%%%%%%%%%%%%%%%%%%%%%%%%%%%%%%%%%%%%%%%%%%%%%%%%%%%%%%%%%%%%%%%%%%%%%%%%%%%
\section{Numerical experiments} \label{numerical_experiments}
%%%%%%%%%%%%%%%%%%%%%%%%%%%%%%%%%%%%%%%%%%%%%%%%%%%%%%%%%%%%%%%%%%%%%%%%%%%%%%%%
All experiments in this article are performed on the Eiger cluster of 
the Swiss National Supercomputing Centre. Each node is 
equipped with 2 x AMD 7742 CPU@2.2 GHz processors.
The calculations are performed using up to 64 MPI processes, 
each consisting of up to 20 OpenMP threads, for a total of up to 1280 cores.
The implementation of the barycentric interpolation has been integrated
into \texttt{Bembel} (the BEM-based engineering library, 
\url{www.bembel.eu}), see \cite{DHK+20}. The latter is also 
used for the assembly of the boundary integral formulation.
%%%%%%%%%%%%%%%%%%%%%%%%%%%%%%%%%%%%%%%%%%%%%%%%%%%%%%%%%%%%%%%%%%%%%%%%%%%%%%%%
\subsection{General setup}
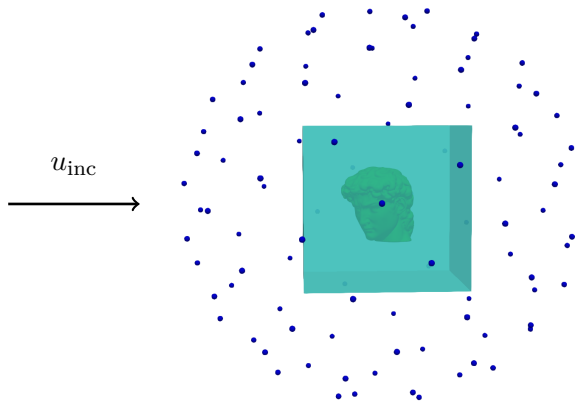
\begin{figure}
	\begin{tikzpicture}[scale=0.25]
		\draw(20,0)node{
			\includegraphics[width=0.5\linewidth,clip,trim=500 200 500 200]{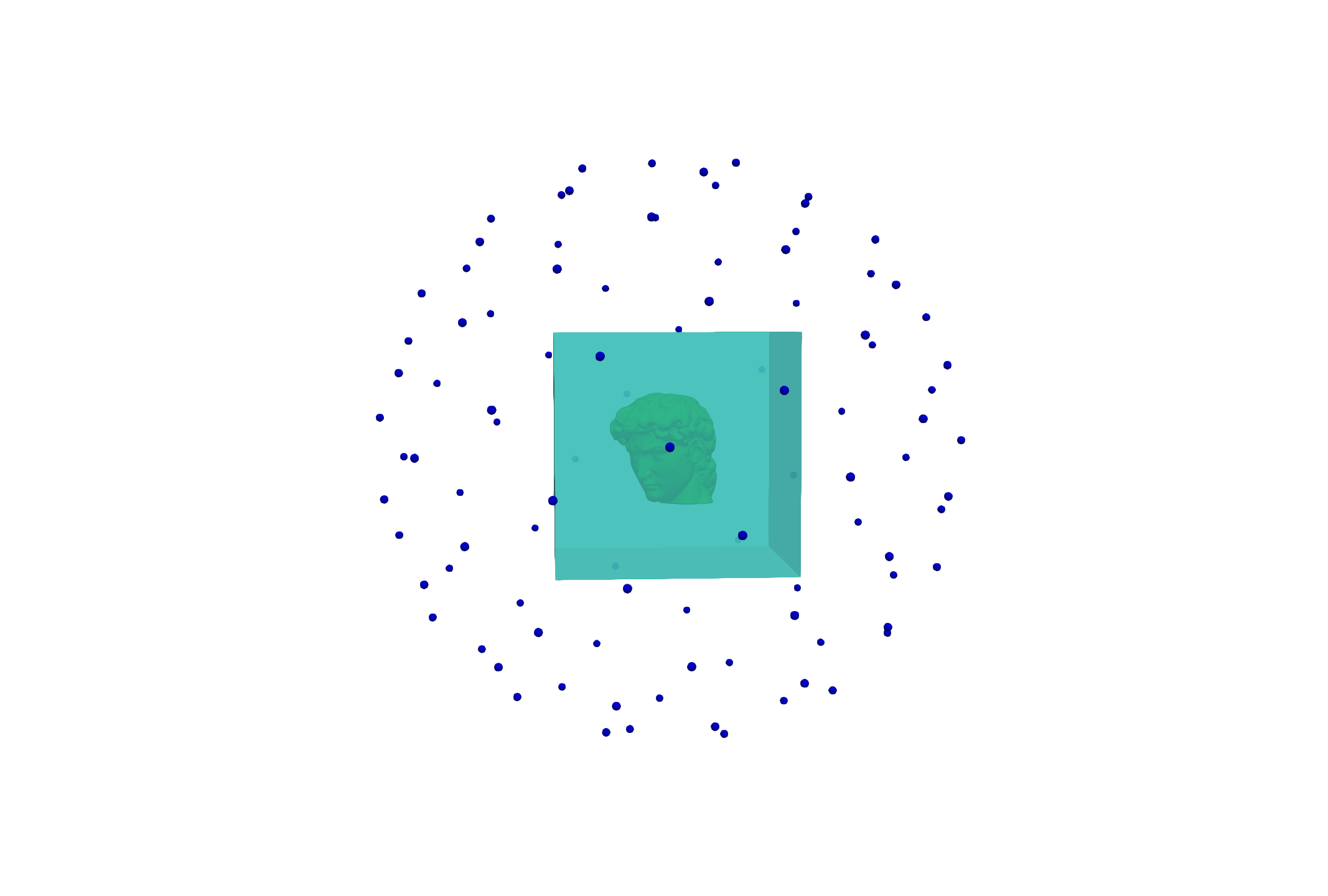}};
		\draw (0,0)node(N0){};
		\draw (8,0)node(N1){};
		\path[->,line width=1pt]  (N0) edge (N1);
		\draw (4,2)node{$u_{\mathrm{inc}}$};
	\end{tikzpicture}
	\caption{Setup for the numerical experiments on the David
 geometry. One hundred blue points are evaluation points. The artificial interface $T$ is the cuboid \([-2,2]^3\). The innermost object inside $T$ is the scatterer. The incident plane wave is considered with wavenumber $\kappa=5$ and direction ${\bs d} = [1,0,0]^\intercal$.}
	\label{fig:setup}
\end{figure}

%%%%%%%%%%%%%%%%%%%%%%%%%%%%%%%%%%%%%%%%%%%%%%%%%%%%%%%%%%%%%%%%%%%%%%%%%%%%%%%%
In the numerical experiments, we consider a torus (bounding box: $[-1,1]\times[-1,1]\times[-1,1]$) 
and Michelangelo's David (bounding box: $[-1,1]\times[-1,1]\times[-0.2,0.2]$) for 
the reference scatterer $D_0$, respectively. Both of domains have a diameter of 2. 
The artificial interface $T$ (24 patches) is given by the cuboid
$[-2,2]\times[-2,2]\times[-2,2]$ in both cases, see Figure~\ref{fig:setup}. The artificial interface undergoes discretization using tensor-product polynomials of degree 8 on each patch. Subsequently, Cauchy data for the artificial interface is derived by conducting 1944 point evaluations on $T$. For the
Helmholtz equation, we consider the incident plane wave with wavenumber 
$\kappa=5$ and direction ${\bs d} = (1,0,0)$. Regarding the random deformation field, 
we consider the deformation field given by 
\(\E[{\bs\chi}](\xref)=\xref\) and
\[
	\Cov[\DefField](\hat{\bs x},\hat{\bs x}')=\begin{bmatrix}
		k_{\frac{3}{2}}(20r)   & 10^{-4} k_{\infty}(4r) & 10^{-4} k_{\infty}(4r) \\
		10^{-4} k_{\infty}(4r) & k_{\frac{3}{2}}(20r)   & 10^{-4} k_{\infty}(4r) \\
		10^{-4} k_{\infty}(4r) & 10^{-4} k_{\infty}(4r) & k_{\frac{3}{2}}(20r)   \\ 
	\end{bmatrix}.
	\]
	
The random field is computed by the pivoted Cholesky decomposition with an accuracy of $10^{-3}$. 
This leads to the parameter dimension \(m=2218, 5563\) for the torus and David, respectively.
Herein, we set the polynomial degree of the barycentric interpolation to 20 for the
torus and to 10 for David. This results in 6516 landmark points for the torus consisting of 16 patches and 18702 landmark points for David consisting of 187 patches. In the deformation model, 
we set the magnitude of the random perturbation $\alpha=0.07$, such that all realisations of 
random scatterers are inside the artificial interface $T$. This results in a maximal possible 
relative displacement of 45\% for the torus and 49\% for David. Three different 
realisations each scatterer are visualised in Figure~\ref{fig:rough_def_torus} and 
Figure~\ref{fig:rough_def_david} for the torus and David, respectively.
	
%=========================================================
\subsection{Boundary element method based on $p$-refinement}
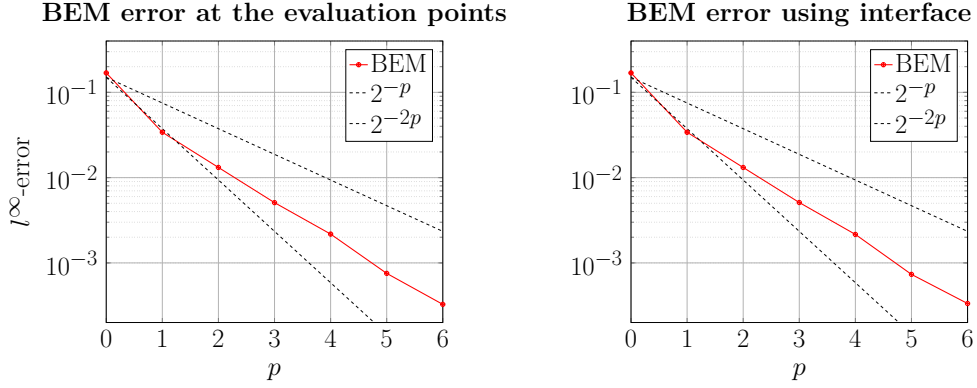
\begin{figure}[htb]
	\begin{center}
		\scalebox{0.4}{
			\pgfplotsset{width=\textwidth}
			\pgfplotsset{minor grid style={dotted,black!30}}
			\pgfplotsset{grid style={solid,black!30}}
			\begin{tikzpicture}
				\Huge
				\begin{semilogyaxis}[ytick distance=10,xtick={0,1,2,3,4,5,6},grid=both, ymin= 2e-4, ymax = 0.4, xmin = 0, xmax = 6,
						legend style={legend pos=north east}, legend cell align={left},%
						xlabel=$p$,
						ylabel={$l^{\infty}$-error},
						title={\textbf{BEM error at the evaluation points}},
					]
					\addplot[line width=1pt,color=red,mark=oplus] table[x index=0,y index=1]{./data/bem/torus/mean_points_1_bem.txt};
					\addplot[dashed] table[x index=0,y expr=0.15*2^-1*\thisrowno{0}]{data/bem/torus/mean_points_1_bem.txt};
					\addplot[dashed] table[x index=0,y expr=0.15*2^-2*\thisrowno{0}]{data/bem/torus/mean_points_1_bem.txt};
					\addlegendentry{BEM};
					\addlegendentry{$2^{-p}$};
					\addlegendentry{$2^{-2p}$};
				\end{semilogyaxis}
			\end{tikzpicture}\qquad\qquad
			\begin{tikzpicture}
				\Huge
				\begin{semilogyaxis}[ytick distance=10,xtick={0,1,2,3,4,5,6},grid=both, ymin= 2e-4, ymax = 0.4, xmin = 0, xmax = 6,
						legend style={legend pos=north east}, legend cell align={left},%
						xlabel=$p$,
						title={\textbf{BEM error using interface}},
					]
					\addplot[line width=1pt,color=red,mark=oplus] table[x index=0,y index=1]{./data/bem/torus/mean_fast_points_1_bem.txt};
					\addplot[dashed] table[x index=0,y expr=0.15*2^-1*\thisrowno{0}]{data/bem/torus/mean_fast_points_1_bem.txt};
					\addplot[dashed] table[x index=0,y expr=0.15*2^-2*\thisrowno{0}]{data/bem/torus/mean_fast_points_1_bem.txt};
					\addlegendentry{BEM};
					\addlegendentry{$2^{-p}$};
					\addlegendentry{$2^{-2p}$};
				\end{semilogyaxis}
			\end{tikzpicture}}
	\end{center}
	\caption{\label{fig:bem_torus}Convergence of boundary element solver on one realisation of the random torus.}
\end{figure}

\begin{figure}
	\begin{center}
		\scalebox{0.4}{
			\pgfplotsset{width=\textwidth}
			\pgfplotsset{minor grid style={dotted,black!30}}
			\pgfplotsset{grid style={solid,black!30}}
			\begin{tikzpicture}
				\Huge
				\begin{semilogyaxis}[ytick distance=10,xtick={0,1,2,3,4,5},grid=both, ymin= 2e-4, ymax = 2e-2, xmin = 0, xmax = 5,
						legend style={legend pos=north east}, legend cell align={left},%
						xlabel=$p$,
						ylabel={$l^{\infty}$-error},
						title={\textbf{BEM error at the evaluation points}},
					]
					\addplot[line width=1pt,color=red,mark=oplus] table[x index=0,y index=1]{./data/bem/david/mean_points_1_bem.txt};
					\addplot[dashed] table[x index=0,y expr=2e-2*2^-1*\thisrowno{0}]{data/bem/david/mean_points_1_bem.txt};
					\addplot[dashed] table[x index=0,y expr=2e-2*2^-2*\thisrowno{0}]{data/bem/david/mean_points_1_bem.txt};
					\addlegendentry{BEM};
					\addlegendentry{$2^{-p}$};
					\addlegendentry{$2^{-2p}$};
				\end{semilogyaxis}
			\end{tikzpicture}\qquad\qquad
			\begin{tikzpicture}
				\Huge
				\begin{semilogyaxis}[ytick distance=10,xtick={0,1,2,3,4,5},grid=both, ymin= 2e-4, ymax = 2e-2, xmin = 0, xmax = 5,
						legend style={legend pos=north east}, legend cell align={left},%
						xlabel=$p$,
						title={\textbf{BEM error using interface}},
					]
					\addplot[line width=1pt,color=red,mark=oplus] table[x index=0,y index=1]{./data/bem/david/mean_fast_points_1_bem.txt};
					\addplot[dashed] table[x index=0,y expr=2e-2*2^-1*\thisrowno{0}]{data/bem/david/mean_fast_points_1_bem.txt};
					\addplot[dashed] table[x index=0,y expr=2e-2*2^-2*\thisrowno{0}]{data/bem/david/mean_fast_points_1_bem.txt};
					\addlegendentry{BEM};
					\addlegendentry{$2^{-p}$};
					\addlegendentry{$2^{-2p}$};
				\end{semilogyaxis}
			\end{tikzpicture}}
	\end{center}
	\caption{\label{fig:bem_david}Convergence of boundary element solver on one realisation of the random David.}
\end{figure}
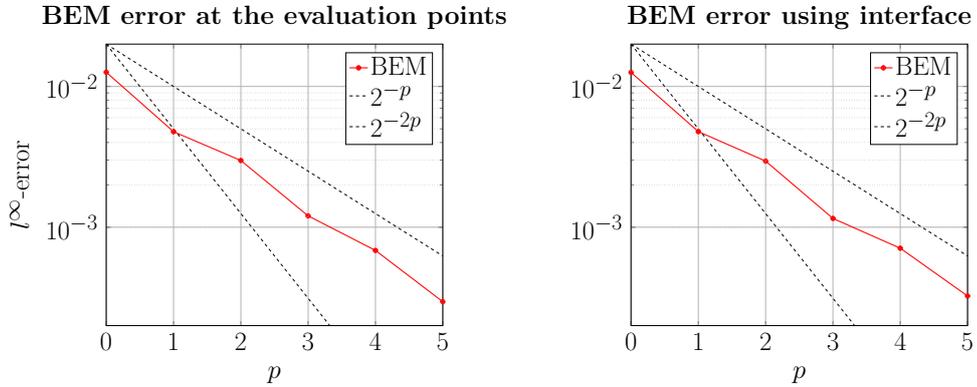

\begin{figure}[htb]
	\begin{center}
		\includegraphics[width=0.3\textwidth,trim=500pt 200pt 550pt 300pt,clip=true]{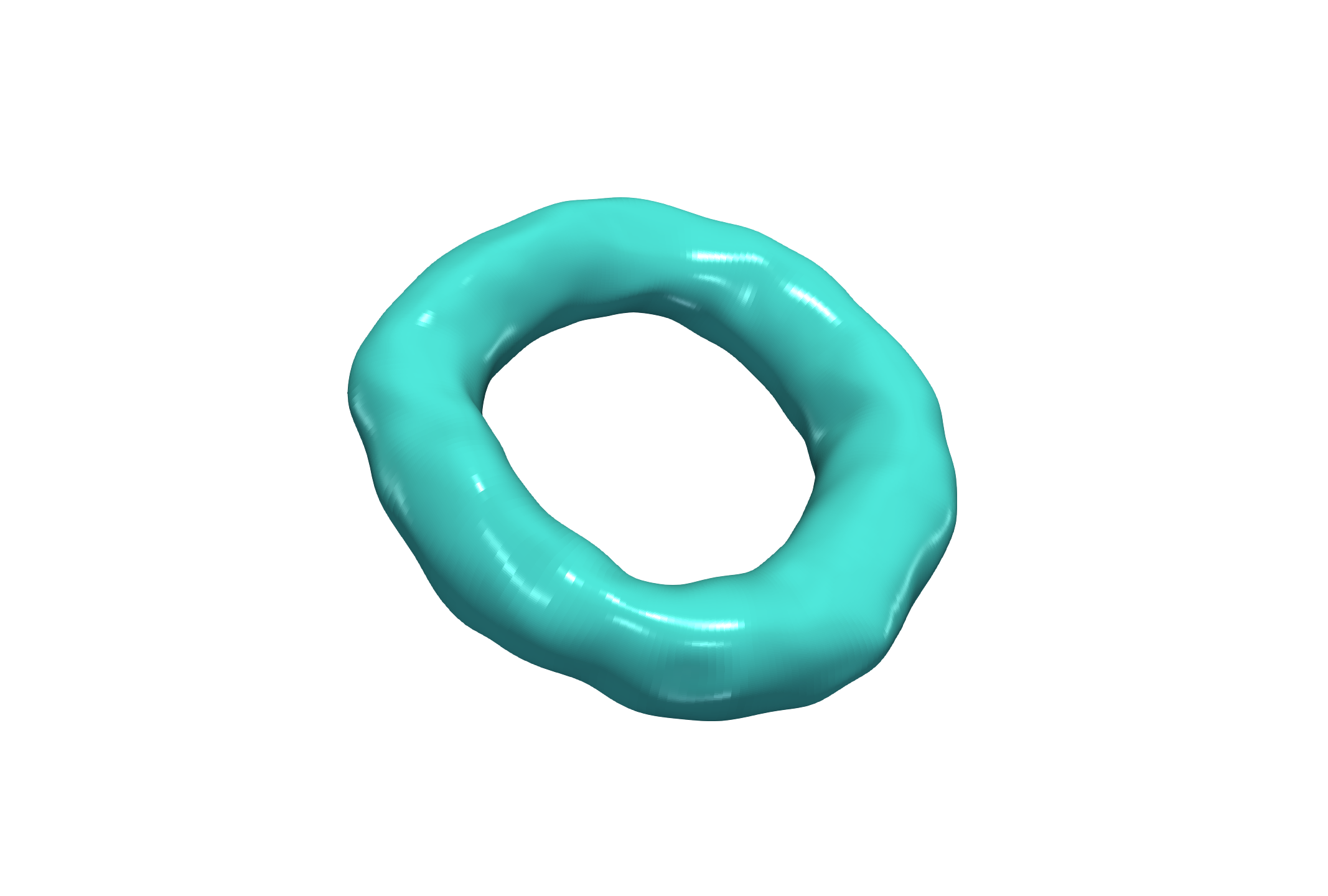}
		\includegraphics[width=0.3\textwidth,trim=500pt 200pt 550pt 300pt,clip=true]{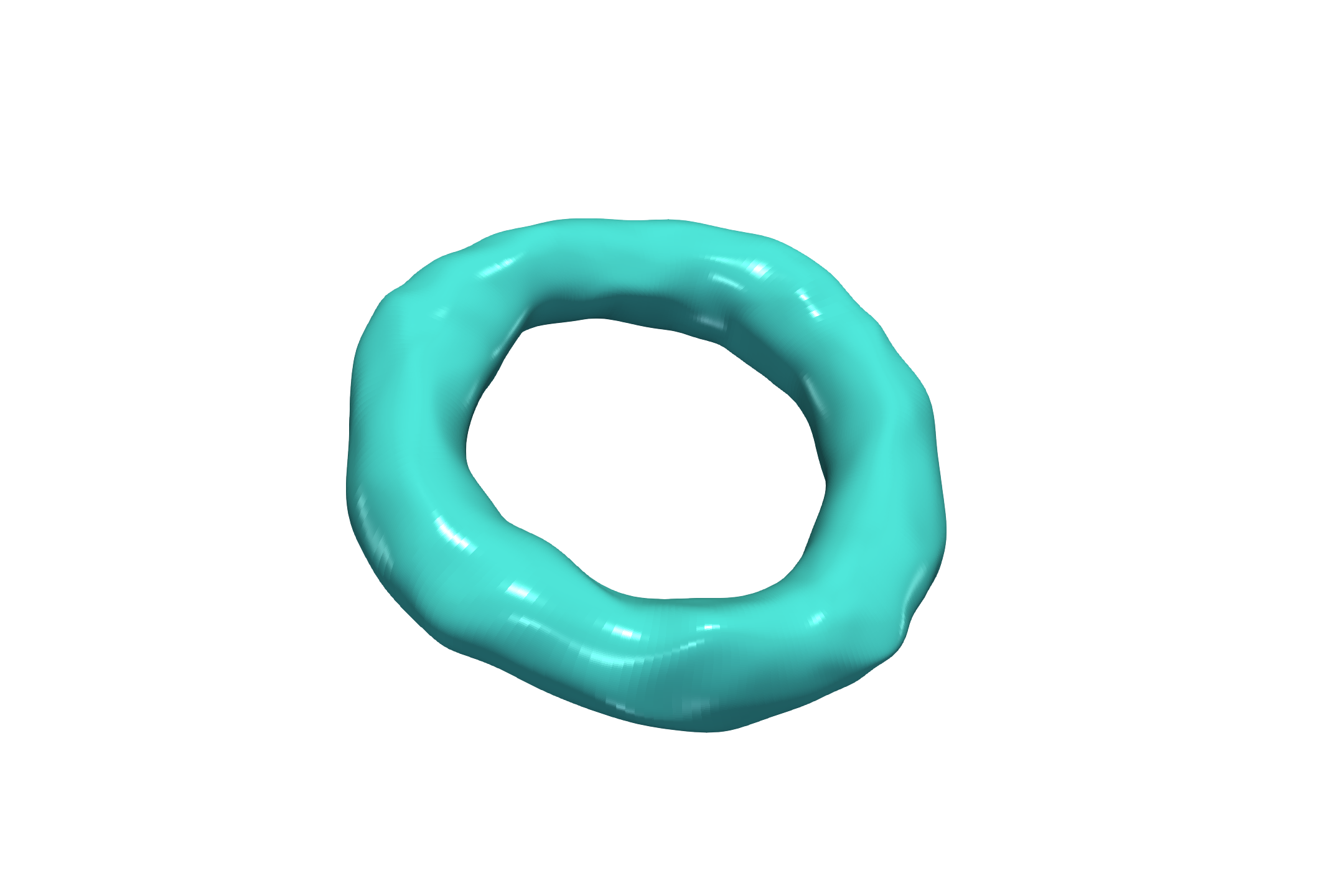}
		\includegraphics[width=0.3\textwidth,trim=500pt 200pt 550pt 300pt,clip=true]{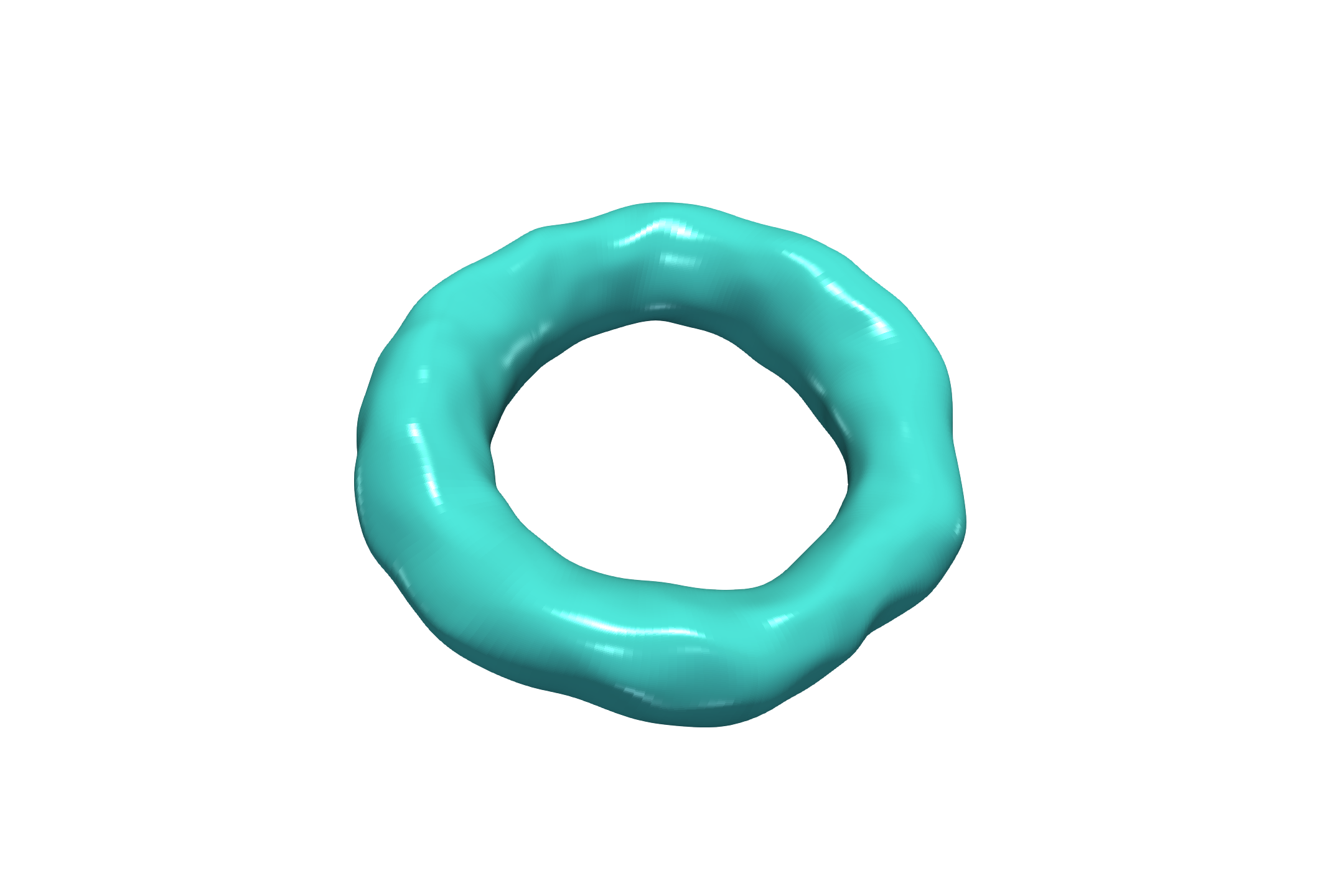}\\
		\includegraphics[width=0.3\textwidth,trim=500pt 0pt 490pt 50pt,clip=true]{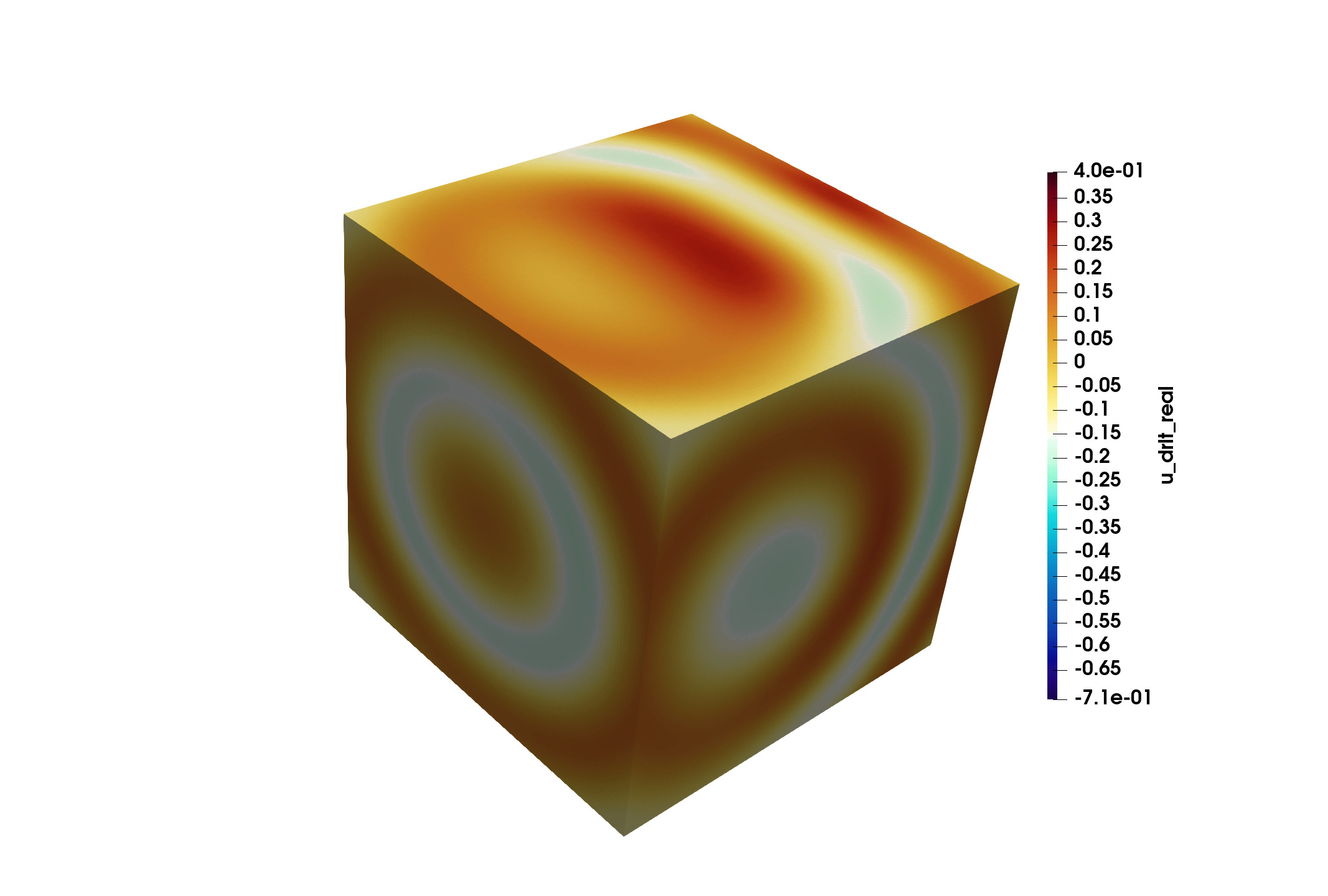}
		\includegraphics[width=0.3\textwidth,trim=500pt 0pt 490pt 50pt,clip=true]{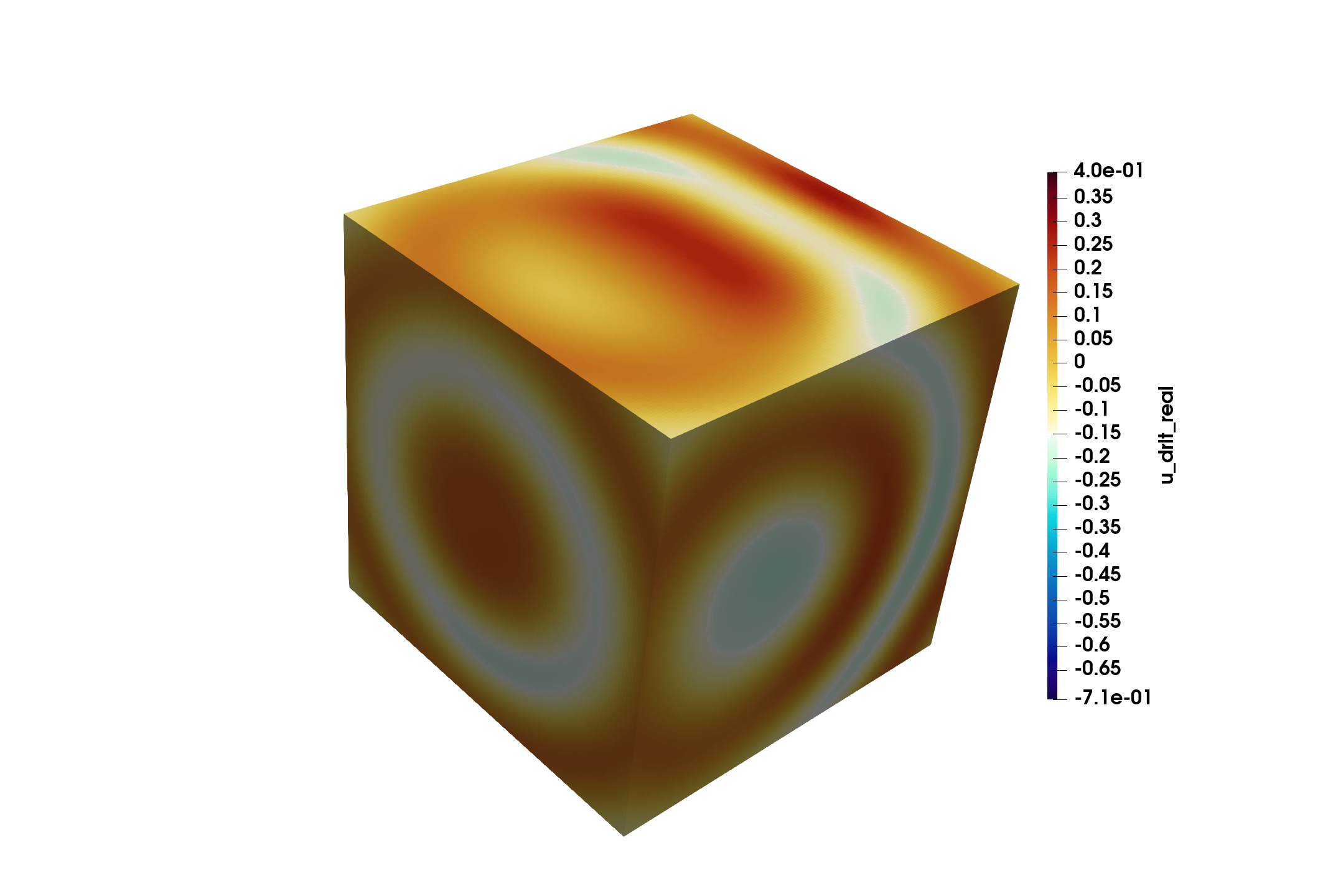}
		\includegraphics[width=0.3\textwidth,trim=500pt 0pt 490pt 50pt,clip=true]{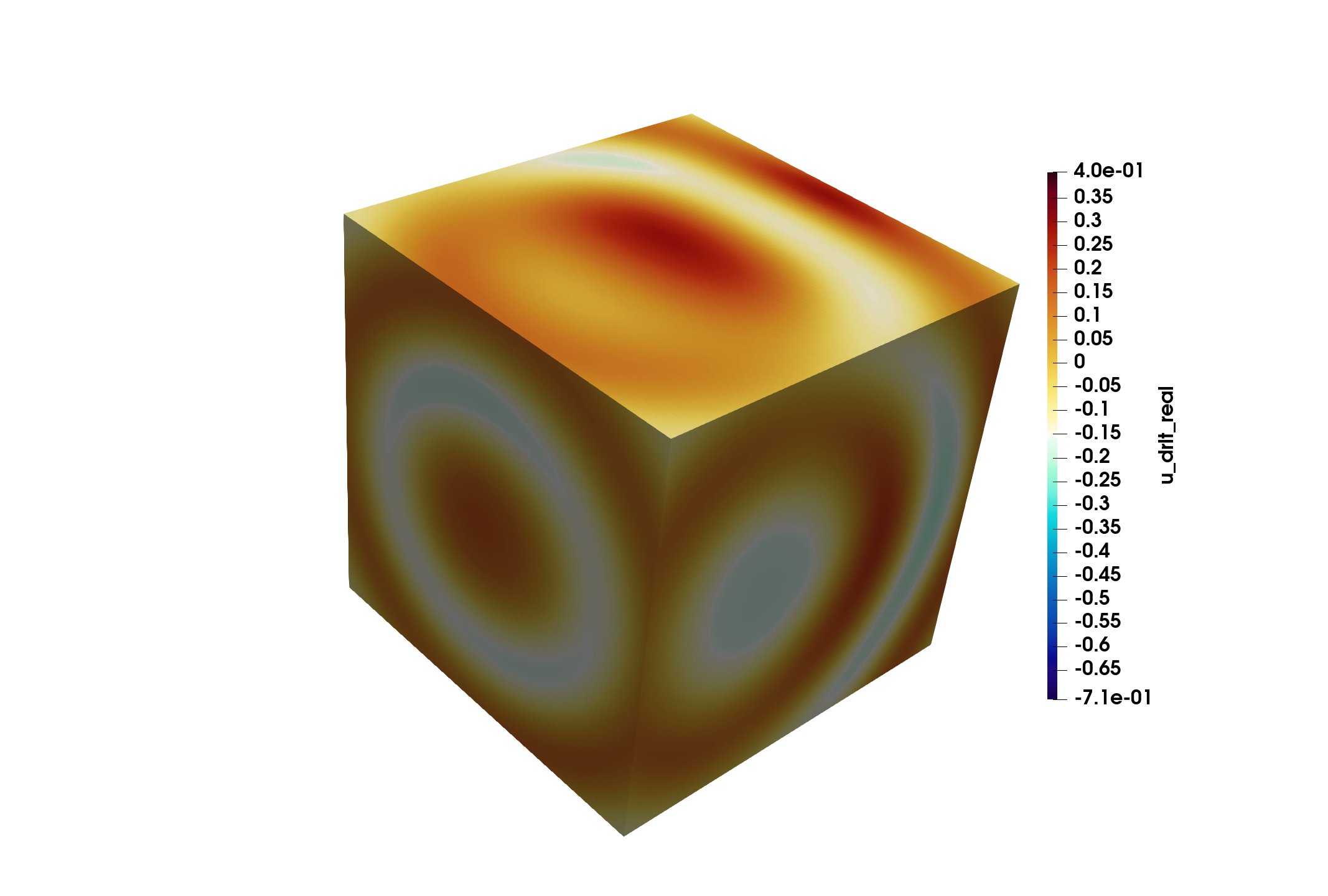}\\
	\end{center}
	\caption{\label{fig:rough_def_torus}Realizations of the scatterer (torus) and the scattered wave at the interface.}
\end{figure}

\begin{figure}[htb]
	\begin{center}
		\includegraphics[width=0.3\textwidth,trim=700pt 200pt 600pt 120pt,clip=true]{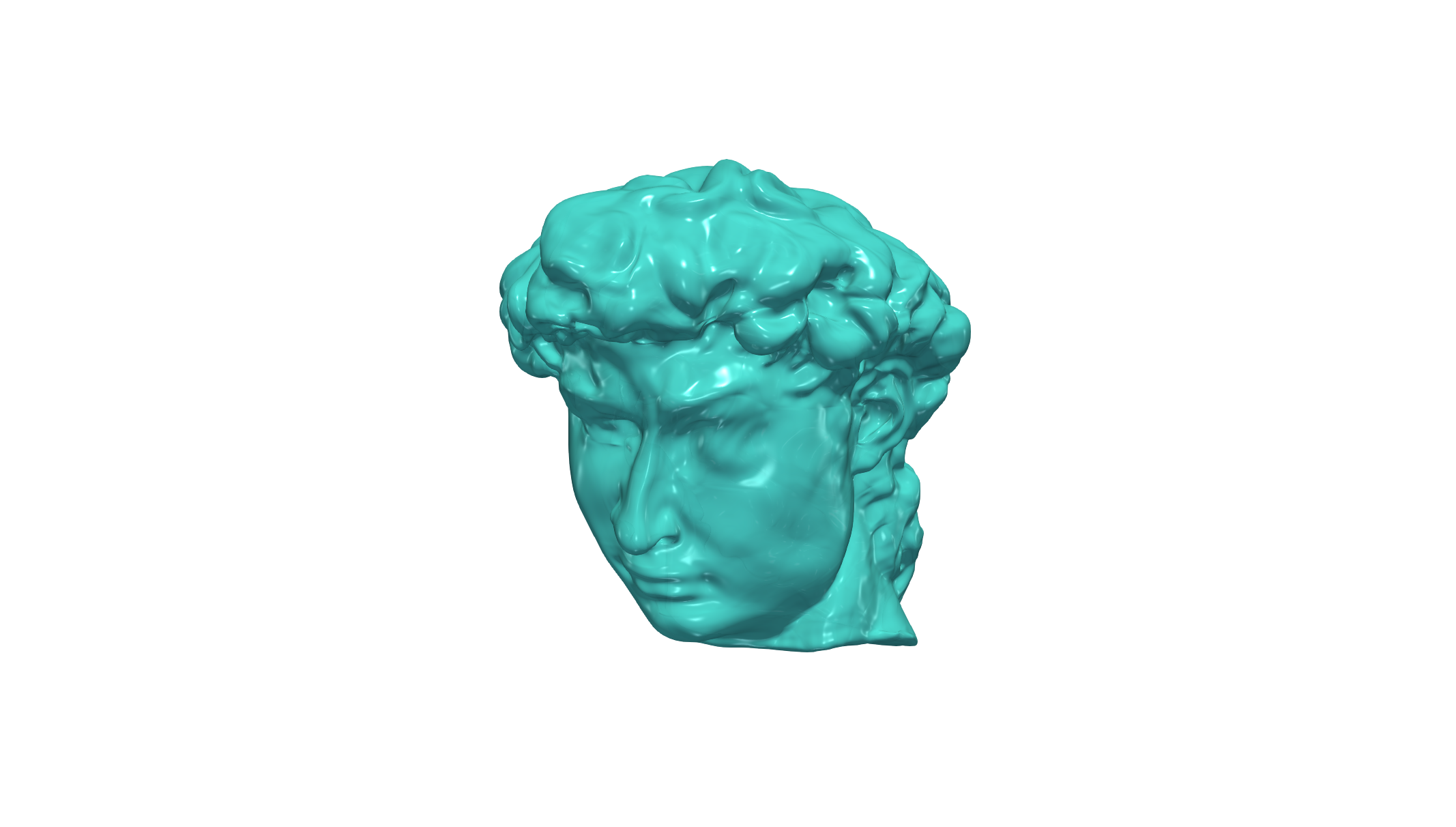}
		\includegraphics[width=0.3\textwidth,trim=700pt 200pt 600pt 120pt,clip=true]{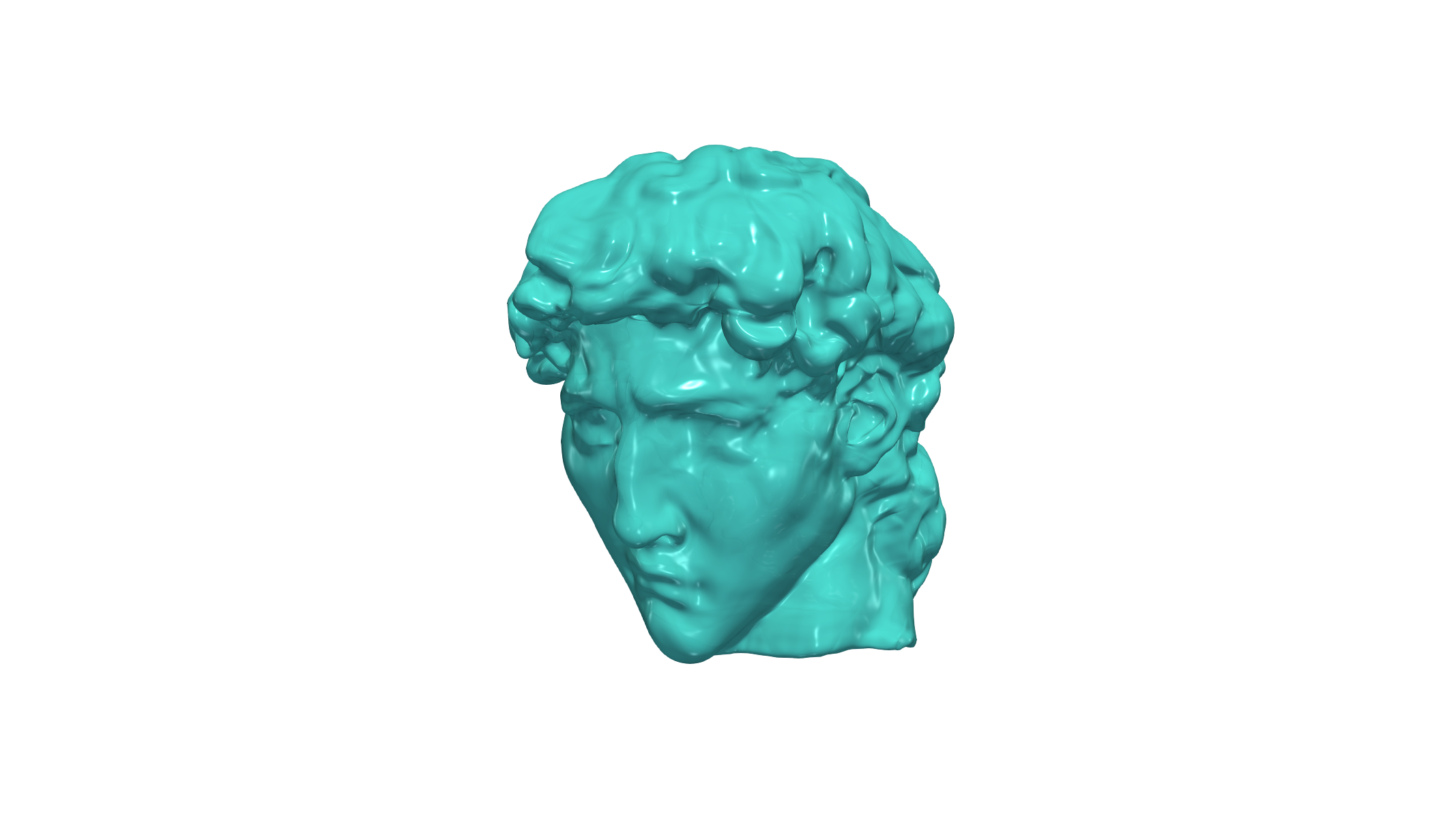}
		\includegraphics[width=0.3\textwidth,trim=700pt 200pt 600pt 120pt,clip=true]{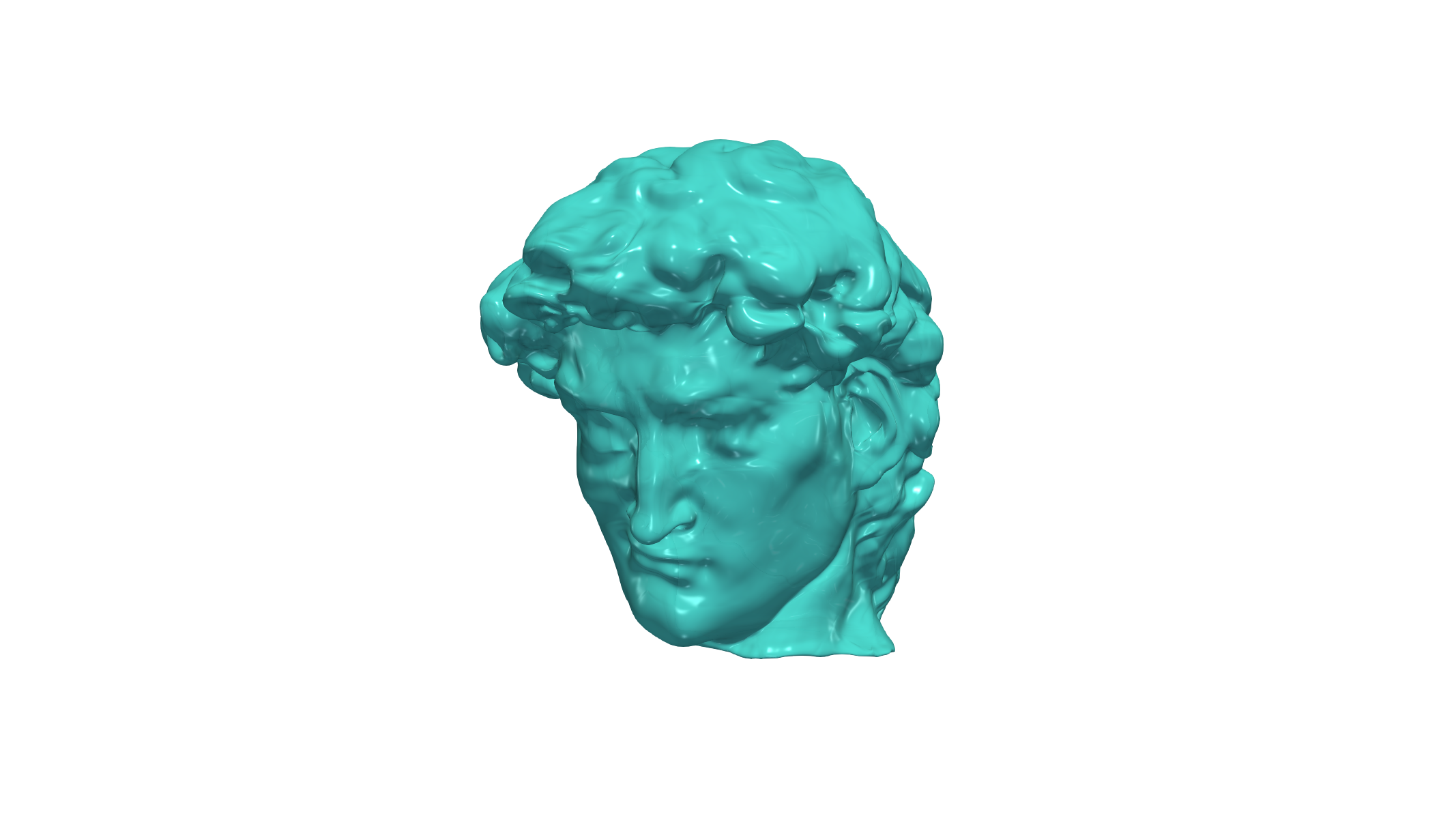}\\
		\includegraphics[width=0.3\textwidth,trim=500pt 50pt 500pt 50pt,clip=true]{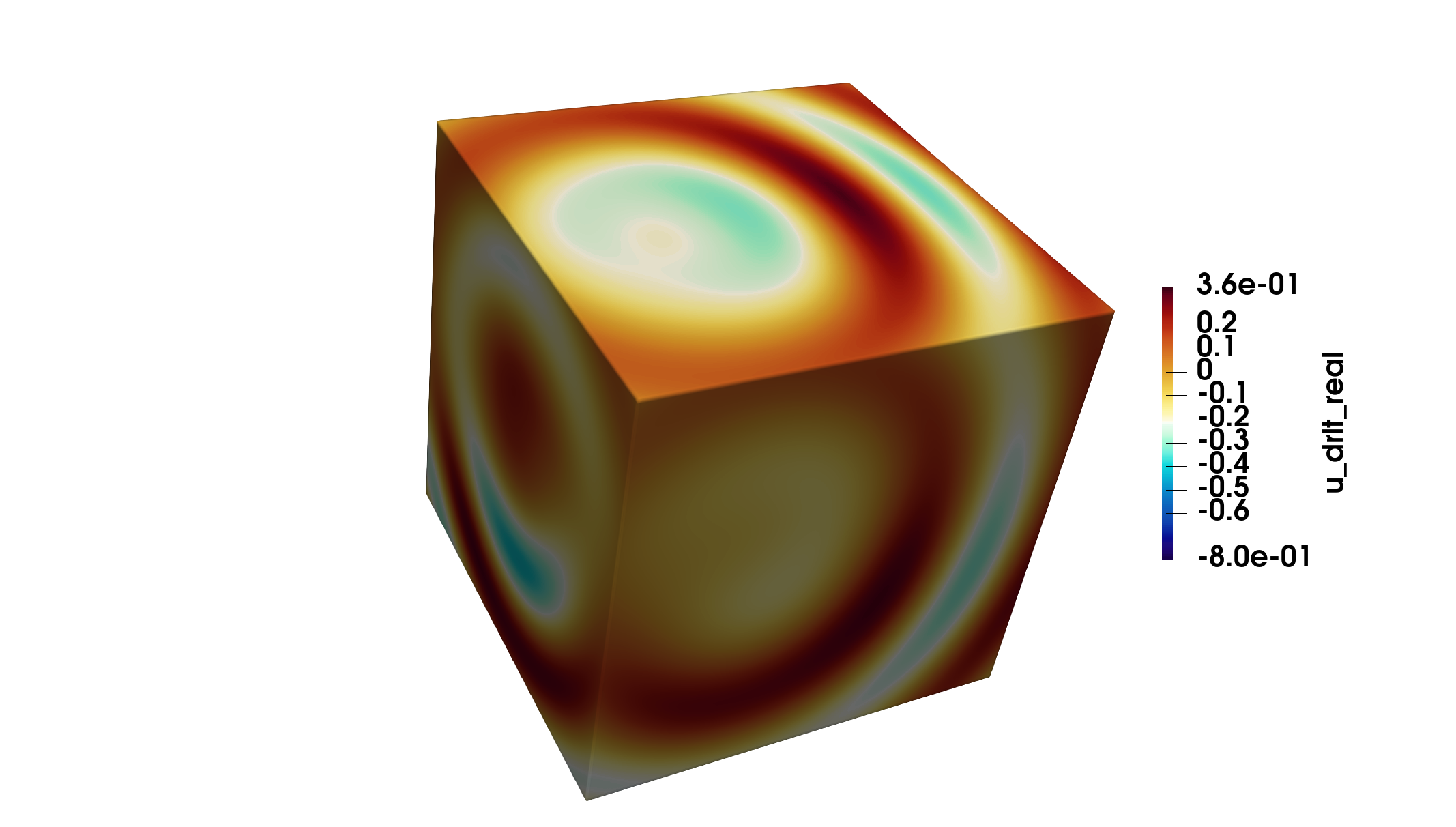}
		\includegraphics[width=0.3\textwidth,trim=500pt 50pt 500pt 50pt,clip=true]{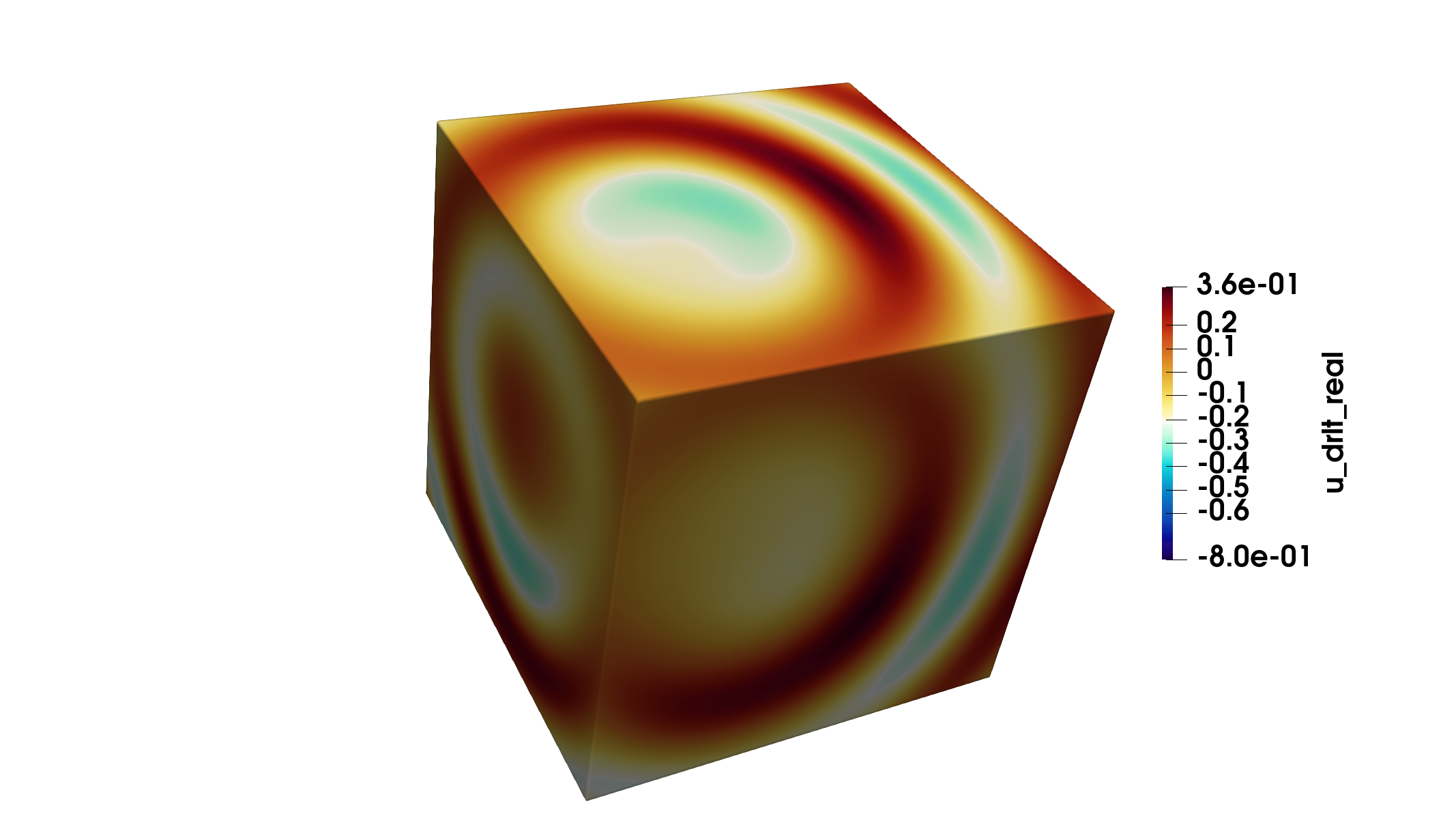}
		\includegraphics[width=0.3\textwidth,trim=500pt 50pt 500pt 50pt,clip=true]{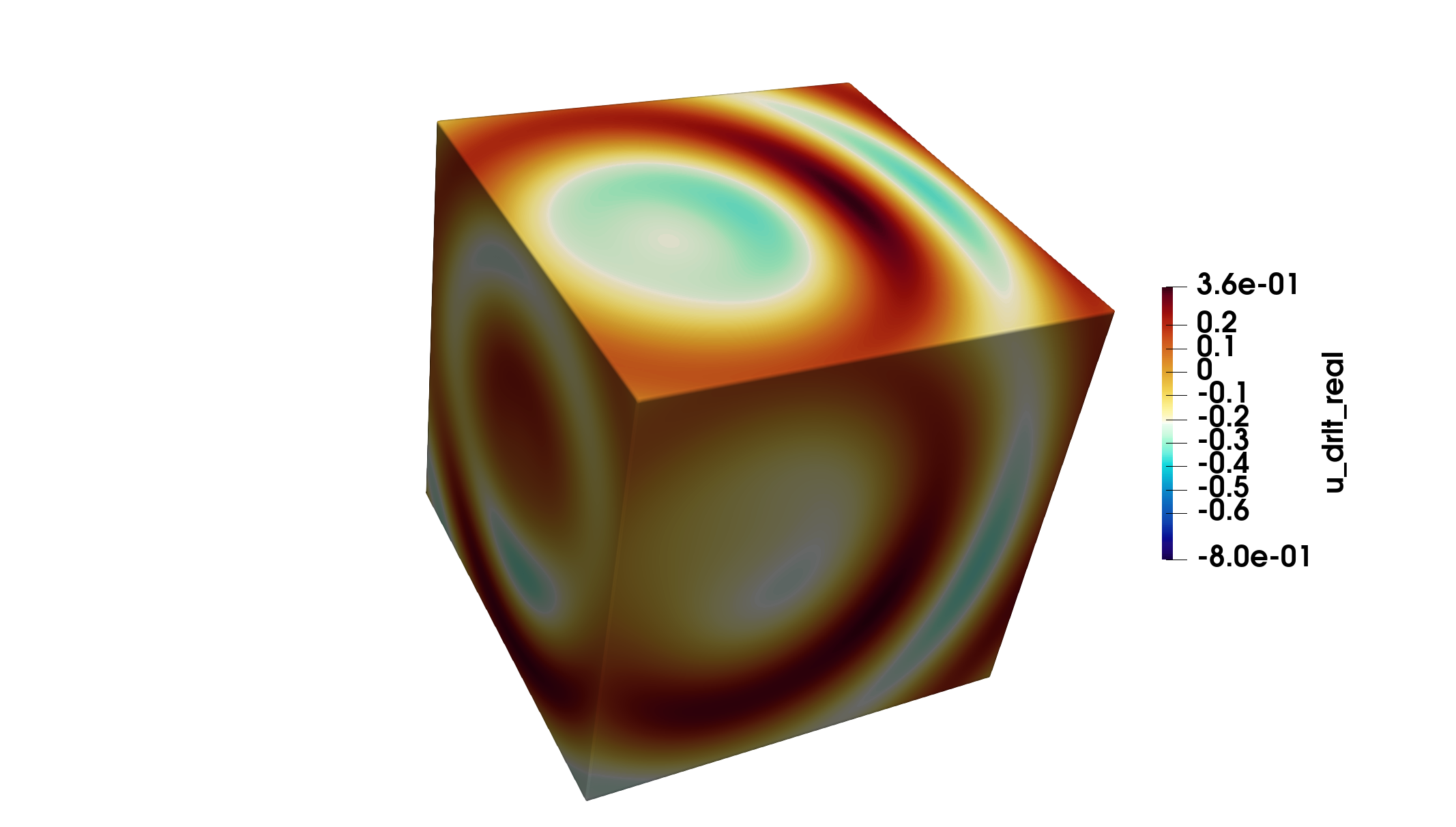}\\
	\end{center}
	\caption{\label{fig:rough_def_david}Realizations of the scatterer (David) and the scattered wave at the interface.}
\end{figure}

%=========================================================
We start by numerically testing the exponential convergence rate of the
maximum potential error evaluated at 100 points in free space with respect
to the polynomial degree $p$. To this end, we consider the first
realization of the random scatterers from Figure~\ref{fig:bem_torus} and Figure~\ref{fig:bem_david}. The evaluation
points are uniformly distributed on a sphere with a radius 5. The number of the
degrees of freedom (DOF) are shown in Table~\ref{tab:tdofs_torus} for the torus
and in Table~\ref{tab:tdofs_david} for David, where we use $p$-refinement at level $\ell=1$. 
When increasing $p$ by one, the degrees of freedom and the computation time in seconds are multiplied by a factor between 2 and 3, which is less than the numbers for
$h$-refinement reported in \cite{dolz2022isogeometric}. The real parts of scattered waves 
by three realizations of each scatterer are depicted in Figure~\ref{fig:rough_def_torus} and Figure~\ref{fig:rough_def_david} for the torus and David, respectively.
	
\begin{table}
	\begin{tabular}{|c|c|c|c|c|c|c|c|}
		\hline
		       & $p=0$    & $p=1$   & $p=2$   & $p=3$   & $p=4$  & $p=5$  & $p=6$  \\\hline
		$t$    & 0.02637s & 0.1139s & 0.3403s & 0.8488s & 1.835s & 3.742s & 7.679s \\\hline
		$DOFs$ & 64       & 144     & 256     & 400     & 576    & 784    & 1024   \\\hline
	\end{tabular}
 	\caption{\label{tab:tdofs_torus}Number of degrees of freedom of the boundary element solver on the torus for different polynomial degrees $p$ and fixed spatial refinement level $\ell=1$.}
\end{table}
	
\begin{table}
	\begin{tabular}{|c|c|c|c|c|c|c|}
		\hline
		       & $p=0$  & $p=1$  & $p=2$ & $p=3$  & $p=4$ & $p=5$  \\\hline
		$t$    & 1.308s & 1.351s & 3.56s & 11.87s & 33.9s & 101.8s \\\hline
		$DOFs$ & 748    & 1683   & 2992  & 4675   & 6732  & 9163   \\\hline
	\end{tabular}
	\caption{\label{tab:tdofs_david}Number of degrees of freedom of the boundary element solver on David for different polynomial degrees $p$ and fixed spatial refinement level $\ell=1$.}
\end{table}
	
%=========================================================
\subsection{$p$-Multilevel Monte Carlo}
\begin{figure}[htb]
	\begin{tikzpicture}
		\huge
		\begin{semilogyaxis}[xmin=0,xmax=2000,ymin=5e-2,ymax=50, ytick distance=10, width=0.7
				\textwidth,height=0.7\textwidth,grid=both,legend style={legend pos=north east, font=\Large}, legend cell align={left},
				xlabel={$k$},
				ylabel={singular value}
			]
			\addplot[mark=*,mark size=.5pt, only marks, color=blue]table[x index=0,y index=1]{./data/evalues_torus.txt};\addlegendentry{Torus};
			\addplot[mark=*,mark size=.5pt, only marks, color=green]table[x index=0,y index=1]{./data/evalues_david.txt};\addlegendentry{David};
			\addplot[red, dashed, line width=1.5pt] table[x index=0,y expr=1000*\thisrowno{0}^-1.25]{./data/evalues_torus.txt};\addlegendentry{$\alpha=1.25$};
		\end{semilogyaxis}
	\end{tikzpicture}
	\caption{\label{fig:eigns} Numerical approximation of the singular values of the 
 Karhunen-Lo\`eve expansion for the model under consideration.}
\end{figure}
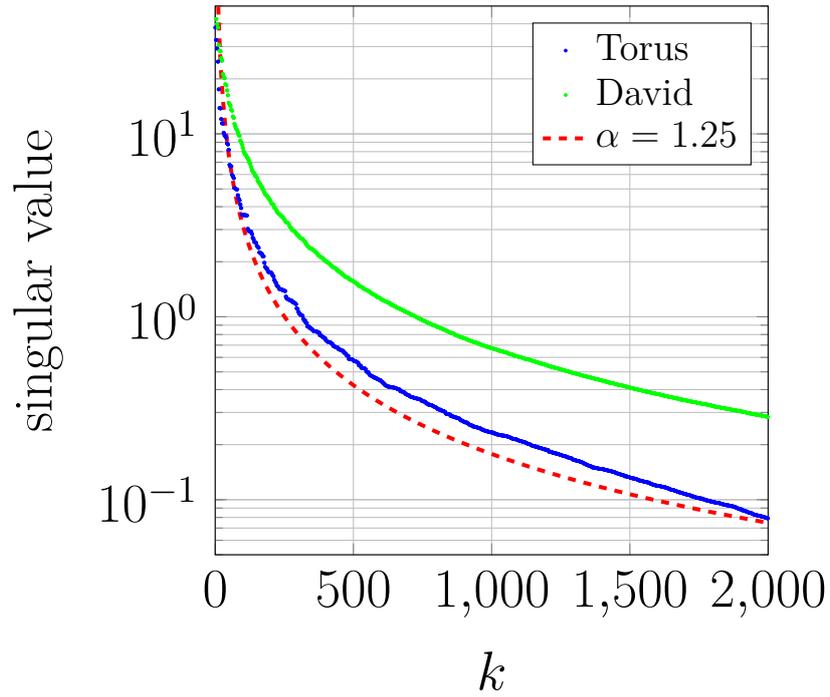

%=========================================================
From \eqref{eq:eig_dec}, the singular values of the chosen Mat\'ern kernel with $\nu=\frac{3}{2}$ asymptotically decay like \(\sigma_k\sim k^{-1.25}\). We numerically 
compute the first 2\,000 largest singular values of the covariance evaluated at the landmark points. The numerical result shows the decay rates are even slower than the theoretical asymptotical rate in the case of David, see Figure~\ref{fig:eigns}. This is caused by the multiplicity \(3\) of each eigenvalue.

\begin{comment}
\begin{figure}[htb]
	\begin{tikzpicture}
		\huge
		\begin{semilogyaxis}[ytick distance=100,xtick={1,2,4,6,8,10,12},grid=both, ymin= 2e-4, ymax = 1.0, xmin = 1, xmax = 12,
				legend style={legend pos=south west, font=\Large}, legend cell align={left},%
				xlabel={sample level},
				ylabel={relative error}
			]
			\addplot[line width=1pt,color=blue,mark=oplus] table[x index=0,y index=2]{./data/qmc_conv_test.txt};\addlegendentry{QMC};
			\addplot[line width=1pt,color=cyan,mark=oplus] table[x index=0,y index=2]{./data/mc_conv_test.txt};\addlegendentry{MC};
			\addplot[dashed] table[x index=0,y expr=.5*2^-0.5*\thisrowno{0}]{data/qmc_conv_test.txt};
			\addlegendentry{$1/\sqrt{N}$};
		\end{semilogyaxis}
	\end{tikzpicture}
	\caption{\label{fig:conv_test} Convergence test of QMC and MC on the high-dimensional integral problem with a slow decaying rate.}
\end{figure}
\end{comment}
	
%%%%%%%%%%%%%%%%%%%%%%%%%%%%%%%%%%%%%%%%%%%%%%%%%%%%%%%%%%%%%%%%%%%%%%%%%%%%%%%%	
We rely on $p$-MLMC to speed up the QoI computations. With increasing levels, 
we decrease the number of samples by a factor of 
\[\hat{\gamma}_p = \sqrt{\frac{\hat{v}_{p-1}}{\hat{v}_p}}\sqrt{\frac{\hat{c}_p}{\hat{c}_{p-1}}}.
\]
Herein, $\hat{v}_p$ and $\hat{c}_p$ are estimated variance and 
time complexity for each polynomial degree $p$, respectively. 
In the experiments, the approximations are obtained using 64 samples. 
The number of samples at the highest polynomial degree $P$ is fixed to 
128. The number of samples increases from $p$ to $p-1$ according to the factor $\hat{\gamma}_p$. 
This leads to the number of samples at each polynomial degree 
level $p$ shown in Table~\ref{tab:samples_torus} and Table~\ref{tab:samples_david} 
for the torus and David, respectively. To benchmark the convergence,
we use the $p$-multilevel Monte Carlo solution with the highest polynomial degree 
$P+1$ as reference. In this study, we set $P$ to 5 and 4 for the torus and David, 
respectively. Notably, the higher number of paths in David's case compensates for its 
lower polynomial degree, in contrast to the torus.
	
	\begin{table}
		\begin{tabular}{|c|c|c|c|c|c|c|c|}
			\hline
			& $p=0$ & $p=1$ & $p=2$ & $p=3$ & $p=4$ & $p=5$ & $p=6$ \\\hline
			reference & 98470 & 47260 & 5969  & 2137  & 1029  & 389   & 128   \\\hline
			$p$-MLMC      & 32390 & 15540 & 1963  & 703   & 338   & 128   & -     \\\hline
		\end{tabular}
  		\caption{\label{tab:samples_torus} Number of samples for the different polynomial degrees $p$ for the $p$-MLMC on the torus. A $p$-MLMC solution with highest degree 6 is used as reference.}
	\end{table}
	
	\begin{table}
		\begin{tabular}{|c|c|c|c|c|c|c|}
			\hline
			& $p=0$  & $p=1$ & $p=2$ & $p=3$ & $p=4$ & $p=5$ \\\hline
			reference & 155500 & 10550 & 2821  & 1491  & 368   & 128   \\\hline
			$p$-MLMC      & 54010  & 3666  & 980   & 518   & 128   & -     \\\hline
		\end{tabular}
		\caption{\label{tab:samples_david}Number of samples for the different polynomial degrees $p$ for the $p$-MLMC on David. A $p$-MLMC solution with degree 5 is used as reference.}
	\end{table}

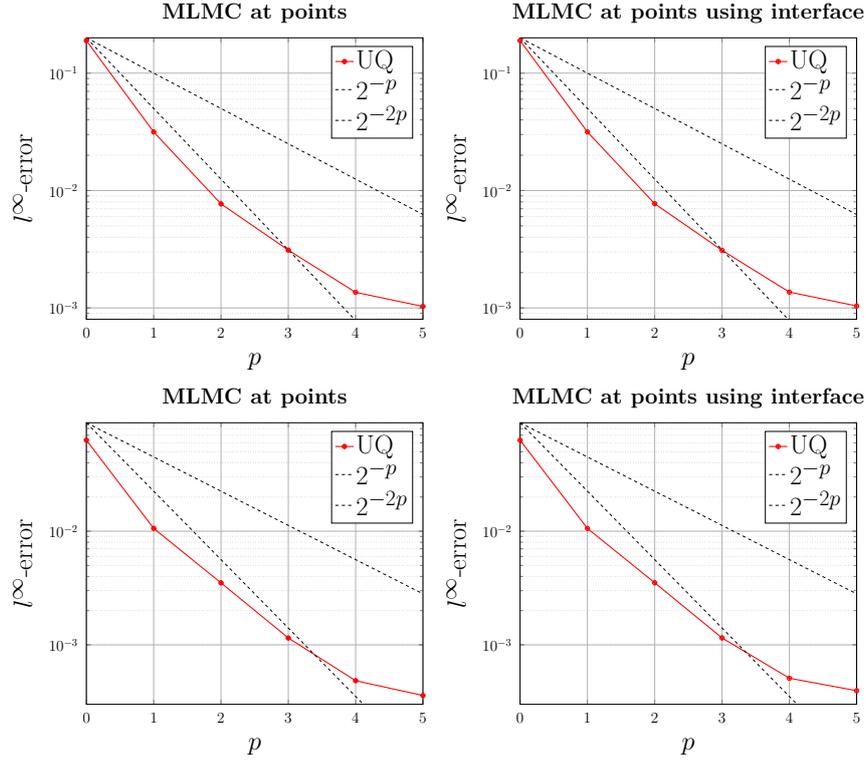
\begin{figure}[htb]
		% Expectation
		\begin{center}
			\scalebox{0.4}{
				\pgfplotsset{width=\textwidth}
				\pgfplotsset{minor grid style={dotted,black!30}}
				\pgfplotsset{grid style={solid,black!30}}
				\begin{tikzpicture}
					\Huge
					\begin{semilogyaxis}[ytick distance=10,xtick={0,1,2,3,4,5},grid=both, ymin= 8e-4, ymax = .2, xmin = 0, xmax = 5,
						legend style={legend pos=north east,font=\Huge}, legend cell align={left},%
						xlabel=$p$,
						ylabel={$l^{\infty}$-error},
						title={\textbf{MLMC at points}},
						xlabel style = {font=\Huge},
						ylabel style = {font=\Huge},
						title style = {font=\huge},
						font=\Large
						]
						\addplot[line width=1pt,color=red,mark=oplus] table[x index=0,y index=1]{./data/uq/torus/mcToref_mean_points_5.txt};
						\addplot[dashed] table[x index=0,y expr=.2*2^-1*\thisrowno{0}]{data/uq/torus/mcToref_mean_interface_5.txt};
						\addplot[dashed] table[x index=0,y expr=.2*2^-2*\thisrowno{0}]{data/uq/torus/mcToref_mean_interface_5.txt};
						\addlegendentry{UQ};
						\addlegendentry{$2^{-p}$};
						\addlegendentry{$2^{-2p}$};
						
					\end{semilogyaxis}
				\end{tikzpicture} 
				
				\begin{tikzpicture}
					\Huge
					\begin{semilogyaxis}[ytick distance=10,xtick={0,1,2,3,4,5},grid=both, ymin= 8e-4, ymax = .2, xmin = 0, xmax = 5,
						legend style={legend pos=north east,font=\Huge}, legend cell align={left},%
						xlabel=$p$,
						ylabel={$l^{\infty}$-error},
						title={\textbf{MLMC at points using interface}},
						xlabel style = {font=\Huge},
						ylabel style = {font=\Huge},
						title style = {font=\huge},
						font=\Large
						]
						\addplot[line width=1pt,color=red,mark=oplus] table[x index=0,y index=1]{./data/uq/torus/mcToref_mean_fast_points_5.txt};
						\addplot[dashed] table[x index=0,y expr=.2*2^-1*\thisrowno{0}]{data/uq/torus/mcToref_mean_interface_5.txt};
						\addplot[dashed] table[x index=0,y expr=.2*2^-2*\thisrowno{0}]{data/uq/torus/mcToref_mean_interface_5.txt};
						\addlegendentry{UQ};
						\addlegendentry{$2^{-p}$};
						\addlegendentry{$2^{-2p}$};
						
					\end{semilogyaxis}
				\end{tikzpicture}
			}
			\scalebox{.4}{
				\pgfplotsset{width=\textwidth}
				\pgfplotsset{minor grid style={dotted,black!30}}
				\pgfplotsset{grid style={solid,black!30}}
				
				\begin{tikzpicture}
					\Huge
					\begin{semilogyaxis}[ytick distance=10,xtick={0,1,2,3,4,5}, grid=both, ymin= 3e-4, ymax = .09, xmin = 0, xmax = 5,
						legend style={legend pos=north east,font=\Huge}, legend cell align={left},%
						xlabel=$p$,
						ylabel={$l^{\infty}$-error},
						title={\textbf{MLMC at points}},
						xlabel style = {font=\Huge},
						ylabel style = {font=\Huge},
						title style = {font=\huge},
						font=\Large
						]
						\addplot[line width=1pt,color=red,mark=oplus] table[x index=0,y index=1]{./data/uq/torus/mcToref_cor_points_5.txt};
						\addplot[dashed] table[x index=0,y expr=.09*2^-1*\thisrowno{0}]{data/uq/torus/mcToref_mean_interface_5.txt};
						\addplot[dashed] table[x index=0,y expr=.09*2^-2*\thisrowno{0}]{data/uq/torus/mcToref_mean_interface_5.txt};
						\addlegendentry{UQ};
						\addlegendentry{$2^{-p}$};
						\addlegendentry{$2^{-2p}$};
						
					\end{semilogyaxis}
				\end{tikzpicture}
				
				\begin{tikzpicture}
					\Huge
					\begin{semilogyaxis}[ytick distance=10,xtick={0,1,2,3,4,5},grid=both, ymin= 3e-4, ymax = .09, xmin = 0, xmax = 5,
						legend style={legend pos=north east,font=\Huge}, legend cell align={left},%
						xlabel=$p$,
						ylabel={$l^{\infty}$-error},
						title={\textbf{MLMC at points using interface}},
						xlabel style = {font=\Huge},
						ylabel style = {font=\Huge},
						title style = {font=\huge},
						font=\Large
						]
						\addplot[line width=1pt,color=red,mark=oplus] table[x index=0,y index=1]{./data/uq/torus/mcToref_cor_fast_points_5.txt};
						\addplot[dashed] table[x index=0,y expr=.09*2^-1*\thisrowno{0}]{data/uq/torus/mcToref_mean_interface_5.txt};
						\addplot[dashed] table[x index=0,y expr=.09*2^-2*\thisrowno{0}]{data/uq/torus/mcToref_mean_interface_5.txt};
						\addlegendentry{UQ};
						\addlegendentry{$2^{-p}$};
						\addlegendentry{$2^{-2p}$};
						
					\end{semilogyaxis}
				\end{tikzpicture}
			}
		\end{center}
		\caption{\label{fig:scatter_err_torus}Convergence of the expectation (top) and correlation (bottom) towards the reference solution for the torus.}
	\end{figure}

\begin{figure}[htb]
	% Expectation
	\begin{center}
		\scalebox{0.4}{
			\pgfplotsset{width=\textwidth}
			\pgfplotsset{minor grid style={dotted,black!30}}
			\pgfplotsset{grid style={solid,black!30}}
			\begin{tikzpicture}
				\Huge
				\begin{semilogyaxis}[ytick distance=10,xtick={0,1,2,3,4},grid=both, ymin= 5e-4, ymax = 2e-2, xmin = 0, xmax = 4,
						legend style={legend pos=north east,font=\Huge}, legend cell align={left},%
						xlabel=$p$,
						ylabel={$l^{\infty}$-error},
						title={\textbf{MLMC at points}},
						xlabel style = {font=\Huge},
						ylabel style = {font=\Huge},
						title style = {font=\huge},
						font=\Large
					]
					\addplot[line width=1pt,color=red,mark=oplus] table[x index=0,y index=1]{./data/uq/david/mcToref_mean_points_4.txt};
					\addplot[dashed] table[x index=0,y expr=2e-2*2^-1*\thisrowno{0}]{data/uq/david/mcToref_mean_interface_4.txt};
					\addplot[dashed] table[x index=0,y expr=2e-2*2^-2*\thisrowno{0}]{data/uq/david/mcToref_mean_interface_4.txt};
					\addlegendentry{UQ};
					\addlegendentry{$2^{-p}$};
					\addlegendentry{$2^{-2p}$};
											
				\end{semilogyaxis}
			\end{tikzpicture} 
							
			\begin{tikzpicture}
				\Huge
				\begin{semilogyaxis}[ytick distance=10,xtick={0,1,2,3,4,5},grid=both, ymin= 5e-4, ymax = 2e-2, xmin = 0, xmax = 4,
						legend style={legend pos=north east,font=\Huge}, legend cell align={left},%
						xlabel=$p$,
						ylabel={$l^{\infty}$-error},
						title={\textbf{MLMC at points using interface}},
						xlabel style = {font=\Huge},
						ylabel style = {font=\Huge},
						title style = {font=\huge},
						font=\Large
					]
					\addplot[line width=1pt,color=red,mark=oplus] table[x index=0,y index=1]{./data/uq/david/mcToref_mean_fast_points_4.txt};
					\addplot[dashed] table[x index=0,y expr=2e-2*2^-1*\thisrowno{0}]{data/uq/david/mcToref_mean_interface_4.txt};
					\addplot[dashed] table[x index=0,y expr=2e-2*2^-2*\thisrowno{0}]{data/uq/david/mcToref_mean_interface_4.txt};
					\addlegendentry{UQ};
					\addlegendentry{$2^{-p}$};
					\addlegendentry{$2^{-2p}$};
											
				\end{semilogyaxis}
			\end{tikzpicture}
		}
		\scalebox{.4}{
			\pgfplotsset{width=\textwidth}
			\pgfplotsset{minor grid style={dotted,black!30}}
			\pgfplotsset{grid style={solid,black!30}}
			\begin{tikzpicture}
				\Huge
				\begin{semilogyaxis}[ytick distance=10,xtick={0,1,2,3,4,5},grid=both, ymin= 2e-4, ymax = 2e-2, xmin = 0, xmax = 4,
						legend style={legend pos=north east,font=\Huge}, legend cell align={left},%
						xlabel=$p$,
						ylabel={$l^{\infty}$-error},
						title={\textbf{MLMC at points}},
						xlabel style = {font=\Huge},
						ylabel style = {font=\Huge},
						title style = {font=\huge},
						font=\Large
					]
					\addplot[line width=1pt,color=red,mark=oplus] table[x index=0,y index=1]{./data/uq/david/mcToref_cor_points_4.txt};
					\addplot[dashed] table[x index=0,y expr=2e-2*2^-1*\thisrowno{0}]{data/uq/david/mcToref_mean_interface_4.txt};
					\addplot[dashed] table[x index=0,y expr=2e-2*2^-2*\thisrowno{0}]{data/uq/david/mcToref_mean_interface_4.txt};
					\addlegendentry{UQ};
					\addlegendentry{$2^{-p}$};
					\addlegendentry{$2^{-2p}$};
											
				\end{semilogyaxis}
			\end{tikzpicture}
							
			\begin{tikzpicture}
				\Huge
				\begin{semilogyaxis}[ytick distance=10,xtick={0,1,2,3,4,5},grid=both, ymin= 2e-4, ymax = 2e-2, xmin = 0, xmax = 4,
						legend style={legend pos=north east,font=\Huge}, legend cell align={left},%
						xlabel=$p$,
						ylabel={$l^{\infty}$-error},
						title={\textbf{MLMC at points using interface}},
						xlabel style = {font=\Huge},
						ylabel style = {font=\Huge},
						title style = {font=\huge},
						font=\Large
					]
					\addplot[line width=1pt,color=red,mark=oplus] table[x index=0,y index=1]{./data/uq/david/mcToref_cor_fast_points_4.txt};
					\addplot[dashed] table[x index=0,y expr=2e-2*2^-1*\thisrowno{0}]{data/uq/david/mcToref_mean_interface_4.txt};
					\addplot[dashed] table[x index=0,y expr=2e-2*2^-2*\thisrowno{0}]{data/uq/david/mcToref_mean_interface_4.txt};
					\addlegendentry{UQ};
					\addlegendentry{$2^{-p}$};
					\addlegendentry{$2^{-2p}$};
											
				\end{semilogyaxis}
			\end{tikzpicture}
		}
	\end{center}
	\caption{\label{fig:scatter_err_david}Convergence of the expectation (top) and correlation (bottom) towards the reference solution for David.}
\end{figure}
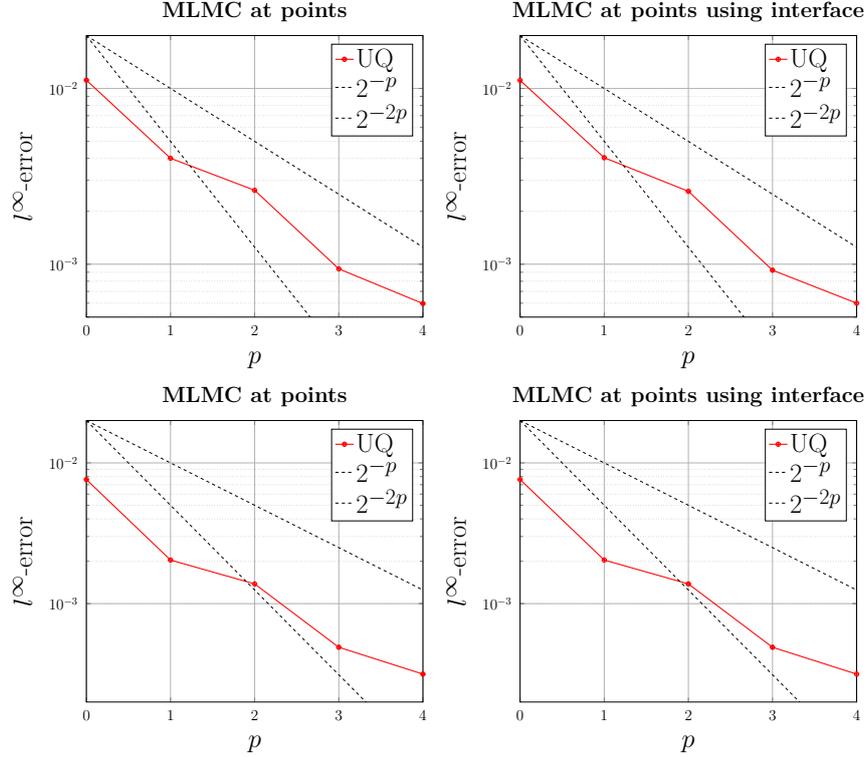

From Figure~\ref{fig:scatter_err_torus} (for the torus) and 
Figure~\ref{fig:scatter_err_david} (for David), the maximum error of 
the expectation and correlation at 100 evaluation points is reduced by a factor approximately between 2 and 4 when increasing $p$ by $1$ (applicable for most values of $p$) in both cases. We report both, the error for the direct evaluation and the evaluation
using the artificial interface. Both errors are comparable. The top row
of each figure corresponds to the expectation and the bottom row to the
correlation. In all cases, the overall convergence rate is almost the same as for the boundary element solver alone. And the maximum error of the expectation at the highest level $P$ approaches a similar scale (around $10^{-3}$) as the error for the boundary element solver alone.
	
%%%%%%%%%%%%%%%%%%%%%%%%%%%%%%%%%%%%%%%%%%%%%%%%%%%%%%%%%%%%%%%%%%%%%%%%%%%%%%%%
\section{Conclusion} \label{conclusion}
In this article, we considered the time-harmonic acoustic scattering problem 
from objects subject to large and rough random deviations. In this setting, 
neither QMC nor sparse grid approaches for quadrature based computation of 
statistical moments are applicable. Perturbation approaches are not suitable 
due to their restriction to rather small deviations. To overcome the prohibitive 
computational cost the combination of Monte Carlo sampling and highly refined 
meshes, we have considered a $p$-multilevel Monte Carlo approach. To this end, 
we have described how we may generate conforming multi-patch NURBS surfaces based 
on complex geometries given by triangular meshes. We have provided
an algorithm by means of barycentric interpolation that allows to
efficiently model even rough deformation fields with low Sobolev regularity of 
the covariance function and a slowly decaying
Karhunen-Lo\`eve expansion. Thanks to the suggested
boundary integral approach, we do not require nested meshes to evaluate the
approximations for different polynomial degrees.
Our numerical studies on complex geometries clearly demonstrate the advantages 
of this approach, in particular the ability of handling large and rough random 
perturbations and the capability of approximating the scattered wave in free 
space.
%%%%%%%%%%%%%%%%%%%%%%%%%%%%%%%%%%%%%%%%%%%%%%%%%%%%%%%%%%%%%%%%%%%%%%%%%%%%%%%%
\bibliography{literature}

\begin{thebibliography}{10}

\bibitem{AJZ20}
R.~Aylwin, C.~Jerez-Hanckes, C.~Schwab, and J.~Zech.
\newblock Domain uncertainty quantification in computational electromagnetics.
\newblock {\em SIAM/ASA Journal on Uncertainty Quantification}, 8(1):301--341,
  2020.

\bibitem{berrut2004barycentric}
J.-P. Berrut and L.~N. Trefethen.
\newblock Barycentric lagrange interpolation.
\newblock {\em SIAM review}, 46(3):501--517, 2004.

\bibitem{ttk19}
T.~Bin~Masood, J.~Budin, M.~Falk, G.~Favelier, C.~Garth, C.~Gueunet,
  P.~Guillou, L.~Hofmann, P.~Hristov, A.~Kamakshidasan, C.~Kappe, P.~Klacansky,
  P.~Laurin, J.~Levine, J.~Lukasczyk, D.~Sakurai, M.~Soler, P.~Steneteg,
  J.~Tierny, W.~Usher, J.~Vidal, and M.~Wozniak.
\newblock {An overview of the Topology ToolKit}.
\newblock In {\em {IEEE Workshop on Topological Data Analysis and
  Visualization}}, 2019.

\bibitem{BRFL+20}
P.~Blondeel, P.~Robbe, C.~Van~hoorickx, S.~François, G.~Lombaert, and
  S.~Vandewalle.
\newblock p-refined multilevel quasi-monte carlo for galerkin finite element
  methods with applications in civil engineering.
\newblock {\em Algorithms}, 13(5):110, 2020.

\bibitem{BN2014}
F.~Bonizzoni and F.~Nobile.
\newblock Perturbation analysis for the darcy problem with log-normal
  permeability.
\newblock {\em SIAM/ASA Journal on Uncertainty Quantification}, 2(1):223--244,
  Jan. 2014.

\bibitem{BG2004a}
H.-J. Bungartz and M.~Griebel.
\newblock Sparse grids.
\newblock {\em Acta Numerica}, 13:147--269, May 2004.

\bibitem{Caf98}
R.~Caflisch.
\newblock Monte {C}arlo and quasi-{M}onte {C}arlo methods.
\newblock {\em Acta Numerica}, 7:1--49, 1998.

\bibitem{Caf1998}
R.~E. Caflisch.
\newblock Monte {{Carlo}} and quasi-{{Monte Carlo}} methods.
\newblock {\em Acta Numerica}, 7:1--49, Jan. 1998.

\bibitem{campen2017partitioning}
M.~Campen.
\newblock Partitioning surfaces into quad patches.
\newblock In {\em Eurographics (Tutorials)}, 2017.

\bibitem{CNT2016}
J.~E. {Castrill{\'o}n-Cand{\'a}s}, F.~Nobile, and R.~F. Tempone.
\newblock Analytic regularity and collocation approximation for elliptic
  {{PDEs}} with random domain deformations.
\newblock {\em Computers \& Mathematics with Applications}, 71(6):1173--1197,
  2016.

\bibitem{CNT2021}
J.~E. {Castrill{\'o}n-Cand{\'a}s}, F.~Nobile, and R.~F. Tempone.
\newblock A hybrid collocation-perturbation approach for {{PDEs}} with random
  domains.
\newblock {\em Advances in Computational Mathematics}, 47(3):40, June 2021.

\bibitem{CK2}
D.~Colton and R.~Kress.
\newblock {\em Inverse Acoustic and Electromagnetic Scattering}.
\newblock Springer, Berlin-Heidelberg-New York, 2nd edition, 1997.

\bibitem{darmon2020acoustic}
M.~Darmon, V.~Dorval, and F.~Baqu{\'e}.
\newblock Acoustic scattering models from rough surfaces: A brief review and
  recent advances.
\newblock {\em Applied Sciences}, 10(22):8305, 2020.

\bibitem{Dol2020}
J.~D{\"o}lz.
\newblock A higher order perturbation approach for electromagnetic scattering
  problems on random domains.
\newblock {\em SIAM/ASA Journal on Uncertainty Quantification}, 8(2):748--774,
  Jan. 2020.

\bibitem{DH2018}
J.~D{\"o}lz and H.~Harbrecht.
\newblock Hierarchical matrix approximation for the uncertainty quantification
  of potentials on random domains.
\newblock {\em Journal of Computational Physics}, 371:506--527, Oct. 2018.

\bibitem{dolz2022isogeometric}
J.~D{\"o}lz, H.~Harbrecht, C.~Jerez-Hanckes, and M.~Multerer.
\newblock Isogeometric multilevel quadrature for forward and inverse random
  acoustic scattering.
\newblock {\em Computer Methods in Applied Mechanics and Engineering},
  388:114242, 2022.

\bibitem{DHK+20}
J.~D{\"o}lz, H.~Harbrecht, S.~Kurz, M.~Multerer, S.~Sch{\"o}ps, and F.~Wolf.
\newblock Bembel: The fast isogeometric boundary element {C++} library for
  {L}aplace, {H}elmholtz, and electric wave equation.
\newblock {\em SoftwareX}, 11:100476, 2020.

\bibitem{DHM23}
J.~D{\"o}lz, H.~Harbrecht, and M.~Multerer.
\newblock Solving acoustic scattering problems by the isogeometric boundary
  element method.
\newblock 2023.

\bibitem{DH2023}
J.~D{\"o}lz and F.~Henr{\'i}quez.
\newblock Parametric shape holomorphy of boundary integral operators with
  applications, May 2023.

\bibitem{dong2006spectral}
S.~Dong, P.-T. Bremer, M.~Garland, V.~Pascucci, and J.~C. Hart.
\newblock Spectral surface quadrangulation.
\newblock In {\em ACM SIGGRAPH 2006 Papers}, pages 1057--1066. 2006.

\bibitem{giles2008multilevel}
M.~B. Giles.
\newblock Multilevel {M}onte {C}arlo path simulation.
\newblock {\em Operations Research}, 56(3):607--617, 2008.

\bibitem{Gil15}
M.~B. Giles.
\newblock Multilevel {M}onte {C}arlo methods.
\newblock {\em Acta Numerica}, 24:259--328, 2015.

\bibitem{graham2015quasi}
I.~G. Graham, F.~Y. Kuo, J.~A. Nichols, R.~Scheichl, C.~Schwab, and I.~H.
  Sloan.
\newblock Quasi-{M}onte {C}arlo finite element methods for elliptic pdes with
  lognormal random coefficients.
\newblock {\em Numerische Mathematik}, 131:329--368, 2015.

\bibitem{HHPS18}
A.-L. Haji-Ali, H.~Harbrecht, M.~D. Peters, and M.~Siebenmorgen.
\newblock Novel results for the anisotropic sparse grid quadrature.
\newblock {\em Journal of Complexity}, 47:62--85, 2018.

\bibitem{HPS}
H.~Harbrecht, M.~Peters, and R.~Schneider.
\newblock On the low-rank approximation by the pivoted {C}holesky
  decomposition.
\newblock {\em Applied Numerical Mathematics}, 62(4):428--440, 2012.

\bibitem{HPS12}
H.~Harbrecht, M.~Peters, and M.~Siebenmorgen.
\newblock On multilevel quadrature for elliptic stochastic partial differential
  equations.
\newblock In J.~Garcke and M.~Griebel, editors, {\em Sparse Grids and
  Applications}, volume~88 of {\em Lecture Notes in Computational Science and
  Engineering}, pages 161--179, Berlin-Heidelberg, 2012. Springer.

\bibitem{HPS16}
H.~Harbrecht, M.~Peters, and M.~Siebenmorgen.
\newblock Analysis of the domain mapping method for elliptic diffusion problems
  on random domains.
\newblock {\em Numerische Mathematik}, 134(4):823--856, 2016.

\bibitem{HSS2008a}
H.~Harbrecht, R.~Schneider, and C.~Schwab.
\newblock Sparse second moment analysis for elliptic problems in stochastic
  domains.
\newblock {\em Numerische Mathematik}, 109(3):385--414, 2008.

\bibitem{H2}
S.~Heinrich.
\newblock Multilevel {M}onte {C}arlo methods.
\newblock In {\em Lecture Notes in Large Scale Scientific Computing}, pages
  58--67, London, 2001. Springer.

\bibitem{HSSS2018}
R.~Hiptmair, L.~Scarabosio, C.~Schillings, and {\relax Ch}.~Schwab.
\newblock Large deformation shape uncertainty quantification in acoustic
  scattering.
\newblock {\em Advances in Computational Mathematics}, 44(5):1475--1518, Oct.
  2018.

\bibitem{huang2022isogeometric}
W.~Huang and M.~Multerer.
\newblock Isogeometric analysis of diffusion problems on random surfaces.
\newblock {\em Applied Numerical Mathematics}, 179:50--65, 2022.

\bibitem{JSZ17}
C.~Jerez-Hanckes, C.~Schwab, and J.~Zech.
\newblock Electromagnetic wave scattering by random surfaces: Shape holomorphy.
\newblock {\em Mathematical Models and Methods in Applied Sciences},
  27(12):2229--2259, 2017.

\bibitem{kachoyan1987acoustic}
B.~Kachoyan and C.~Macaskill.
\newblock Acoustic scattering from an arbitrarily rough surface.
\newblock {\em The Journal of the Acoustical Society of America},
  82(5):1720--1726, 1987.

\bibitem{macaskill1988numerical}
C.~Macaskill and B.~Kachoyan.
\newblock Numerical evaluation of the statistics of acoustic scattering from a
  rough surface.
\newblock {\em The Journal of the Acoustical Society of America},
  84(5):1826--1835, 1988.

\bibitem{Mul18}
M.~D. Multerer.
\newblock A note on the domain mapping method with rough diffusion
  coefficients.
\newblock {\em Applied Numerical Mathematics}, 145:283--296, 2019.

\bibitem{salzer1972lagrangian}
H.~E. Salzer.
\newblock Lagrangian interpolation at the chebyshev points xn, $\nu$$\equiv$
  cos ($\nu$$\pi$/n), $\nu$= 0 (1) n; some unnoted advantages.
\newblock {\em The Computer Journal}, 15(2):156--159, 1972.

\bibitem{ttk}
J.~Tierny, G.~Favelier, J.~A. Levine, C.~Gueunet, and M.~Michaux.
\newblock The {T}opology {T}ool{K}it.
\newblock {\em {IEEE Transactions on Visualization and Computer Graphics (Proc.
  of IEEE VIS)}}, 2017.
\newblock \url{https://topology-tool-kit.github.io/}.

\bibitem{twersky1983reflection}
V.~Twersky.
\newblock Reflection and scattering of sound by correlated rough surfaces.
\newblock {\em The Journal of the Acoustical Society of America}, 73(1):85--94,
  1983.

\bibitem{Wan02}
X.~Wang.
\newblock A constructive approach to strong tractability using quasi-{M}onte
  {C}arlo algorithms.
\newblock {\em J. Complexity}, 18:683--701, 2002.

\bibitem{xiu}
D.~Xiu and D.~M. Tartakovsky.
\newblock Numerical methods for differential equations in random domains.
\newblock {\em SIAM Journal on Scientific Computing}, 28(3):1167--1185, 2006.

\end{thebibliography}
\bibliographystyle{abbrv}
\end{document}